\documentclass[11pt,a4paper,dvipsnames]{article}
\usepackage[top=25mm, bottom=25mm, right=25mm, left=25mm]{geometry}

\usepackage[T1]{fontenc}
\usepackage[utf8]{inputenc}

\usepackage[english]{babel}      

\usepackage{lmodern}           
\usepackage{microtype}         
\usepackage{mathrsfs}          
\usepackage{amsthm, amsmath, amssymb, mathtools} 
\usepackage{aliascnt}         
\usepackage{bm, bbm}           
\usepackage{stmaryrd}          
\usepackage{accents}
\numberwithin{equation}{section}

\usepackage{cases,xfrac,gensymb} 
\usepackage{tocloft}             
\usepackage{listings}            
\usepackage{lipsum}              
\usepackage{titling}             
\usepackage{proof}               
\let\oldunderbar\underbar
\usepackage{sectsty}             
\let\underbar\oldunderbar
\usepackage{parskip}             
\usepackage{booktabs}            
\usepackage{enumerate}           
\usepackage{autobreak}           
\usepackage[table]{xcolor}       
\usepackage{csquotes}            
\usepackage{adjustbox}

\usepackage[bf,sf]{titlesec}       
\usepackage[labelfont={bf,sf}]{caption} 
\usepackage{subcaption}            

\usepackage[
  backend=biber,   
  maxnames=20,
  giveninits=true,  
  isbn=false,        
]{biblatex}
\AtEveryBibitem{%
  \clearfield{month}%
}

\usepackage{graphicx}            
\usepackage{float}               
\usepackage{tikz}                
\usetikzlibrary{external}

\usepackage{pgfplots}
\pgfplotsset{table/search path={.}}
\pgfplotsset{compat=1.18}
\usepgfplotslibrary{colormaps}
\usepackage{pgfplotstable}
\usepgfplotslibrary{fillbetween}

\let\mathbf\bm
\DeclareMathAlphabet\mathbfcal{OMS}{cmsy}{b}{n} 

\newcommand{\abs}[1]{\left\lvert#1\right\rvert}

\newcommand{\norm}[1]{\lVert#1\rVert}

\newcommand{\inner}[2]{\langle #1, #2 \rangle}

\newcommand{\reals}{\mathbb{R}}

\newcommand{\naturals}{\mathbb{N}_0}

\newcommand{\sym}{\mathbb{S}}

\newcommand{\calH}{\mathcal{H}}

\newcommand{\prox}{{\rm{prox}}}
\newcommand{\kron}{\otimes}

\DeclareMathOperator*{\minimize}{minimize}

\DeclareMathOperator*{\argmin}{argmin}

\DeclareMathOperator*{\gra}{gra}

\DeclareMathOperator*{\zer}{zer}
\DeclareMathOperator*{\dom}{dom}
\DeclareMathOperator*{\Id}{Id}

\newcommand{\quadform}[2]{\mathcal{Q}\p{#1,#2}}

\newcommand{\xmiddle}[1]{\;\middle#1\;}

\newcommand{\bx}{\mathbf{x}}
\newcommand{\bu}{\mathbf{u}}
\newcommand{\by}{\mathbf{y}}
\newcommand{\bz}{\mathbf{z}}
\newcommand{\bfcn}{\mathbf{f}}
\newcommand{\bFcn}{\mathbf{F}}

\newcommand{\IndexOp}{\mathcal{I}_{\textup{op}}}
\newcommand{\IndexFunc}{\mathcal{I}_{\textup{func}}}
\newcommand{\NumFunc}{m_{\textup{func}}}
\newcommand{\NumOp}{m_{\textup{op}}}
\newcommand{\NumEval}{\bar{m}}
\newcommand{\NumEvalOp}{\bar{m}_{\textup{op}}}
\newcommand{\NumEvalFunc}{\bar{m}_{\textup{func}}}

\newcommand\set[1]{\mathord{\left.\{ #1 \} \right. }}
\newcommand\Bigset[1]{\mathord{\left\{ #1 \right\}}}
\newcommand\p[1]{\mathord{( #1 )}}
\newcommand\Bigp[1]{\mathord{\left( #1 \right)}}

\providecommand{\MSC}[1]
{\textbf{\textsf{Mathematics subject classification 2020.}} #1}

\makeatletter
\def\th@plain{
  \thm@headfont{\normalfont\sffamily\bfseries}%
  \itshape
}
\def\th@definition{
  \thm@headfont{\normalfont\sffamily\bfseries}%
  \thm@notefont{\normalfont\sffamily\bfseries}%
}
\newtheoremstyle{myStyle1}
  {0.3cm}
  {0.3cm}
  {\itshape}
  {}
  {}
  {:}
  {.5em}
  {%
    \thmname   {\normalfont\sffamily\bfseries #1 }%
    \thmnumber {\normalfont\sffamily\bfseries #2}%
    \thmnote   {\normalfont\sffamily\bfseries \ (#3)}%
  }

\newtheoremstyle{myStyle2}
  {0.3cm}   
  {0.3cm}   
  {}        
  {}        
  {}
  {:}       
  {.5em}    
  {%
    \thmname   {\normalfont\sffamily\bfseries #1 }%
    \thmnumber {\normalfont\sffamily\bfseries #2}%
    \thmnote   {\normalfont\sffamily\bfseries \ (#3)}%
  }
\makeatother

\theoremstyle{myStyle1}

\newaliascnt{corollary}{theorem}

\aliascntresetthe{corollary}
\newaliascnt{lemma}{theorem}

\aliascntresetthe{lemma}

\theoremstyle{myStyle2}
\newaliascnt{remark}{theorem}
\newtheorem{remark}[remark]{Remark}
\aliascntresetthe{remark}
\newaliascnt{example}{theorem}

\aliascntresetthe{example}

\providecommand{\keywords}[1]{\textbf{\textsf{Keywords.}} #1}

\definecolor{Red}{rgb}{1, 0, 0}
\definecolor{BlueIMT}{RGB}{0,42,72}
\definecolor{ForestGreen}{rgb}{0.13, 0.55, 0.13}
\definecolor{Gray}{rgb}{0.66, 0.66, 0.66}
\definecolor{Green}{rgb}{0.0, 0.5, 0.0}
\definecolor{MidnightBlue}{rgb}{0.098, 0.098, 0.439}
\definecolor{Orange}{rgb}{0.93, 0.53, 0.18}

\definecolor{color1}{RGB}{68,119,170}   
\definecolor{color2}{RGB}{102,204,238}  
\definecolor{color3}{RGB}{34,136,51}    
\definecolor{color4}{RGB}{17,119,51}    
\definecolor{color5}{RGB}{204,187,68}   
\definecolor{color6}{RGB}{221,204,119}  
\definecolor{color7}{RGB}{204,102,119}  
\definecolor{color8}{RGB}{136,34,85}    
\definecolor{color9}{RGB}{170,51,119}   
\definecolor{color10}{RGB}{102,102,102} 
\definecolor{color11}{RGB}{50,50,50}

\usepackage{hyperref}
\hypersetup{%
    hidelinks,
    hypertexnames = true,
    plainpages    = false,
    colorlinks = true,
    urlcolor   = MidnightBlue,
    linkcolor  = MidnightBlue,
    citecolor  = ForestGreen,
}
\usepackage[capitalize,nameinlink]{cleveref}

\theoremstyle{myStyle1}
\newaliascnt{assumption}{theorem}

\aliascntresetthe{assumption}
\newaliascnt{proposition}{theorem}

\aliascntresetthe{proposition}
\newaliascnt{definition}{theorem}
\newtheorem{definition}[definition]{Definition}
\aliascntresetthe{definition}

\crefname{theorem}{Theorem}{Theorems}
\crefname{corollary}{Corollary}{Corollaries}
\crefname{lemma}{Lemma}{Lemmas}
\crefname{remark}{Remark}{Remarks}
\crefname{example}{Example}{Examples}
\crefname{proposition}{Proposition}{Propositions}
\crefname{assumption}{Assumption}{Assumptions}
\crefname{definition}{Definition}{Definitions}

\lstdefinelanguage{Julia}%
  {morekeywords={abstract,break,case,catch,const,continue,do,else,elseif,%
      end,export,false,for,function,immutable,import,importall,if,in,%
      macro,module,otherwise,quote,return,switch,true,try,type,typealias,%
      using,while, begin},%
   sensitive=true,%
   morecomment=[l]\#,%
   morecomment=[n]{\#=}{=\#},%
   morestring=[s]{"}{"},%
   morestring=[m]{'}{'},%
}[keywords,comments,strings]%

\lstdefinestyle{juliastyle}{language=Julia}

\title{\sffamily\bfseries The AutoLyap software suite for computer-assisted Lyapunov analyses of first-order methods \\[2ex]}

\author{
    Manu Upadhyaya$^{1,2}$ \and
    Shuvomoy Das Gupta$^{3}$ \and
    Adrien B. Taylor$^{1}$ \and
    Sebastian Banert$^{4}$ \and
    Pontus Giselsson$^{2}$
    \\[2ex]
}

\date{%
    $^{1}$Inria, D.I. ENS, CNRS, PSL Research University, Paris, France \\
    \texttt{\{\href{mailto:manu.upadhyaya@inria.fr}{manu.upadhyaya},\href{mailto:adrien.taylor@inria.fr}{adrien.taylor}\}@inria.fr} \\[2ex]
    $^{2}$Department of Automatic Control, Lund University, Lund, Sweden \\
    \texttt{\href{mailto:pontus.giselsson@control.lth.se}{pontus.giselsson@control.lth.se}} \\[2ex]
    $^{3}$Department of Computational Applied Mathematics \& Operations Research,\\
    Rice University, Houston, USA \\
    \texttt{\href{sd158@rice.edu}{sd158@rice.edu}} \\[2ex]
    $^{4}$Center for Industrial Mathematics, University of Bremen, Bremen, Germany \\
    \texttt{\href{mailto:banert@uni-bremen.de}{banert@uni-bremen.de}}
}

\begin{document}

\maketitle

\begin{abstract}
    \noindent
We introduce AutoLyap, a software suite that assists with Lyapunov analyses of a wide class of first-order methods for structured optimization and inclusion problems.
Lyapunov analyses are structured proof patterns, with historical roots in the study of dynamical systems, commonly used to establish convergence results for first-order methods.
Building on previous work, the core idea behind AutoLyap is to recast the verification of the existence of a Lyapunov analysis as a semidefinite program~(SDP), which can then be solved numerically using standard SDP solvers.
Users of the package specify (i)~the class of optimization or inclusion problems, (ii)~the first-order method in question, and (iii)~the type of Lyapunov analysis they wish to test.
Once these inputs are provided, AutoLyap handles the SDP modeling and proceeds to solve the SDP numerically.
We use the package to numerically verify and extend several convergence results.
AutoLyap is currently available in \texttt{Python} and \texttt{Julia}.

\end{abstract}


\MSC{
65K05, 
65K10, 
90C22, 
90C25, 
93D30 
}


\section{Introduction}

Lyapunov analyses have become fundamental for establishing structured and compact convergence guarantees for first-order optimization and operator-splitting methods.
Despite their theoretical strength, traditional Lyapunov analyses often demand intricate, manual derivations, limiting their accessibility.
To address these challenges, this paper introduces \texttt{AutoLyap} and \texttt{AutoLyap.jl}, native \texttt{Python} and \texttt{Julia} packages, respectively, that automate the numerical search for Lyapunov analyses via semidefinite programming.
AutoLyap streamlines the verification and derivation of convergence properties by numerically solving the associated SDP formulations. In particular, AutoLyap enables researchers and practitioners to quickly and reliably obtain convergence results for a wide range of structured optimization and inclusion problems.

The documentation of \texttt{AutoLyap} (\texttt{Python} version) can be found at
\begin{center}
\url{https://autolyap.github.io}
\end{center}
The version considered in this paper is \texttt{AutoLyap} v$0$.$2$.$0$.
The \texttt{Python} package supports two interchangeable optimization backends: the proprietary \texttt{MOSEK Fusion API for Python} and the open-source modeling layer \texttt{CVXPY}.
The Fusion backend provides a direct interface to \texttt{MOSEK}~\cite{mosek2025pythonfusionapimanual}, which offers free academic licenses, while the \texttt{CVXPY} backend provides a solver-agnostic interface that can access a broad range of modern SDP solvers such as \texttt{Clarabel}~\cite{clarabel2024interiorpointsolverconicprograms} (open-source), \texttt{MOSEK}~\cite{mosek2025pythonfusionapimanual} (commercial, free for academic use), \texttt{COPT}~\cite{copt2023cardinaloptimizeruserguide} (commercial, free for academic use), \texttt{SCS}~\cite{scs2016conicoptimizationoperatorsplittinghomogeneousselfdualembedding} (open-source), \texttt{COSMO} \cite{garstka2021cosmoconicoperatorsplittingmethod} (open-source), and \texttt{SDPA} \cite{YamashitaEtal2012_SDPAFamily} (open-source).
We also provide \texttt{AutoLyap.jl} v$0$.$1$.$0$ , a native open-source \texttt{Julia} implementation available at
\begin{center}
\url{https://github.com/autolyap/autolyap.jl}
\end{center}
which is built on the solver-agnostic interface \texttt{JuMP}~\cite{lubin2023jumprecentimprovementsmodelinglanguage}. Through \texttt{JuMP}, the Julia implementation provides access to a set of modern open-source and commercial SDP solvers that is largely identical to the one available through the \texttt{CVXPY} backend.


\subsection{Related works}

\subsubsection{Performance estimation problems and related software}

The performance estimation problem (PEP) methodology, first introduced in \cite{drori2014performancefirstorder} and formalized in \cite{taylor2017exactworstcase,taylor2016smoothstronglyconvex}, provides a systematic way to obtain unimprovable (also called tight) performance guarantees for a large class of first-order methods.
The PEP methodology involves finding a worst-case example within a predefined class of problems for the algorithm and performance measure under consideration via an optimization problem referred to as a performance estimation problem.
As in Lyapunov search, such problems can be reformulated as semidefinite programs and solved, up to numerical precision, using existing SDP software.
Closed-form convergence guarantees can then often be extracted from the numerical solution or derived using computer algebra software.
However, the SDPs obtained from vanilla PEPs are typically difficult to handle in practice.
More precisely, if \(K\) denotes the number of iterations of a given first-order algorithm, then the corresponding SDP typically has \(O(k^2)\) variables and \(O(k^2)\) constraints.
While vanilla PEPs have the appealing property of providing tight performance bounds and offering insights into the proof structures of first-order algorithms (see, e.g.,~\cite[\S ``New answers to simple questions'']{taylorconvex} and~\cite{goujaud2023fundamental}), scalability is often a critical issue, particularly when deriving closed-form results.
Another difficulty with the PEP approach is that it naturally focuses on a fixed iteration count \(k\). Thus, the corresponding guarantees do not necessarily generalize beyond this fixed horizon.

These limitations have motivated the development of SDP-based approaches that use Lyapunov techniques to analyze the performance of first-order methods and partially or fully overcome these issues \cite{moucer2023systematicapproachlyapunov,taylor2019stochasticfirstorder,taylor2018lyapunovfunctionsfirst}. This comes at the cost of potentially losing tightness guarantees. In particular, the AutoLyap software suite is based on~\cite{upadhyaya2025automatedtightlyapunov}.
The key idea is to focus on convergence-rate proof patterns based on key Lyapunov inequalities. These can then be reformulated as more tractable but equivalent semidefinite programs, whose number of variables and constraints is either independent of the horizon \(k\) or grows only linearly with \(k\).
Compared to the standard PEP approach, this provides a more tractable framework for deriving closed-form convergence rates and producing proofs that are generally more concise, accessible, and amenable to closed-form analysis.
Lyapunov-based proof patterns thus represent a pragmatic compromise between tractability and tightness for obtaining performance certificates.
They offer a better chance of deriving closed-form convergence rates while still yielding tight guarantees in many settings.
Examples of successful Lyapunov-based analyses in the first-order optimization literature are numerous.
Such examples include the celebrated accelerated gradient method of Nesterov~\cite{nesterov1983fast} and the more recent optimized gradient method (OGM). The latter was first considered in~\cite{drori2014performancefirstorder}, formally obtained in~\cite{kim2015optimizedfirstorder}, shown to be worst-case optimal for minimizing smooth convex functions in~\cite{drori2017exactinformationbased}, and given a tight Lyapunov-based analysis in~\cite[Section 4.3.1]{daspremont2021accelerationmethods}.

Moreover, the original PEP methodology is implemented in the software packages \texttt{PESTO}~\cite{taylor2017performanceestimationtoolbox} (\texttt{MATLAB}) and \texttt{PEPit}~\cite{goujaud2024pepitcomputerassisted} (\texttt{Python}).
While \texttt{PESTO} and \texttt{PEPit} conveniently allow one to verify Lyapunov-based analyses, \texttt{AutoLyap} is specifically designed to search for them.

\subsubsection{Integral quadratic constraints and their relation to AutoLyap}

Another closely related Lyapunov-based approach involves integral quadratic constraints (IQCs), a technique from robust control theory \cite{megretski1997systemanalysisvia}.
IQCs were first adapted for analyzing first-order methods in \cite{lessard2016analysisdesignoptimization} and subsequently extended in various works \cite{fazlyab2018analysisoptimizationalgorithms,jakob2025linearparametervaryingframeworkanalysis,jakob2025onlineconvex,lessard2022analysisoptimizationalgorithms,padmanabhan2025varyingstepsizes,miller2025analysissynthesisswitchedoptimization,scherer2023optimizationalgorithmsynthesis,vanscoy2024speedrobustnesstrade}.
These approaches share a common feature; they represent first-order methods as linear systems interconnected through feedback with nonlinear mappings.
Such representations, widely used in nonlinear systems analysis, offer compact algorithm descriptions.

The methodological developments presented in this work build on both methodologies: the worst-case analysis and tightness guarantees provided by PEP and the compact algorithm representations offered by IQCs.
A first step toward combining these methodologies was taken in \cite{taylor2018lyapunovfunctionsfirst} and further developed in \cite{upadhyaya2025automatedtightlyapunov}.
That line of work provides the foundation for the further formalization, broader applicability, and software accessibility pursued here.

For example, the problem class and algorithm representation considered here generalize those introduced in \cite[Section 2.1]{upadhyaya2025automatedtightlyapunov} and \cite[Section 2.2]{upadhyaya2025automatedtightlyapunov}, respectively.
While the framework presented in \cite{upadhyaya2025automatedtightlyapunov} addresses convex optimization with iteration-independent algorithm parameters, the approach considered here integrates several generalizations previously studied in separate works.
Specifically, we incorporate (i) optimization beyond convex settings, (ii) inclusion problems beyond optimization, (iii) iteration-dependent algorithm parameters, and (iv) non-frugal algorithms, in which basic oracle calls (such as gradient/operator or proximal/resolvent evaluations) may occur multiple times per iteration.

\subsection{Organization}
The rest of the paper is organized as follows.
\Cref{sec:prelim} contains the notation used and mathematical preliminaries. 
In \Cref{sec:problem_and_algorithm}, we detail our modeling approach to structured optimization and inclusion problems, and the algorithms we use to solve them. In \Cref{sec:analysis}, we outline the types of Lyapunov analyses considered. 
Throughout \Cref{sec:problem_and_algorithm,sec:analysis}, we use the \texttt{Python} package \texttt{AutoLyap} as the reference implementation.
In Section~\ref{sec:autolyap_jl}, we discuss the architecture of \texttt{AutoLyap.jl}, explain how its design differs from the \texttt{Python} package \texttt{AutoLyap}, and note that the two packages maintain near-identical syntax for end users.
Finally, \Cref{sec:conclusions} offers concluding remarks and discusses avenues for future research.

\subsection{Notation and preliminaries}\label{sec:prelim}

Let \(\naturals\) denote the set of nonnegative integers, 
\(\mathbb{N}\) the set of positive integers, 
\(\mathbb{Z}\) the set of integers, \(\llbracket n,m \rrbracket = \set{l \in \mathbb{Z} \xmiddle\vert n \leq l \leq m}\) the set of integers between \(n, m \in \mathbb{Z}\cup\set{\pm\infty}\), 
\(\reals\) the set of real numbers, 
\(\reals_+\) the set of nonnegative real numbers, 
\(\reals_{++}\) the set of positive real numbers, 
\(\reals^n\) the set of all \(n\)-tuples of elements of \(\reals\), 
\(\reals^{m\times n}\) the set of real-valued matrices of size \(m\times n\), and 
\(\sym^{n}\) the set of symmetric real-valued matrices of size \(n\times n\).

Throughout this paper, \((\calH,\inner{\cdot}{\cdot})\) will denote a real Hilbert space.
All norms \(\norm{\cdot}\) are the norms induced by the relevant inner products, which will be clear from context.
We denote the identity mapping \(x \mapsto x\) on \(\calH\) by \(\Id\).

\begin{definition}\label{def:func_defs}
Let \(f\colon \calH \to \reals \cup \{\pm\infty\}\), \(L\in\reals_{+}\) and \(\mu,\widetilde{\mu},\mu_{\textup{gd}}\in\reals_{++}\).
The function \(f\) is said to be
    \begin{enumerate}[(i)]
\item \emph{proper} if \(-\infty \notin f\p{\calH}\) and \(\dom f \neq \emptyset\), where the set \(\dom f = \{x\in\calH \allowbreak \mid \allowbreak f(x) \allowbreak < + \infty\}\) is called the \emph{effective domain} of \(f,\)
\item \emph{lower semicontinuous} if $\liminf_{y\to x} f\p{y} \geq f\p{x}$ for each $x\in\calH$,
\item \emph{convex} if \(f\p{\p{1 - \lambda} x + \lambda y} \leq \p{1 - \lambda} f\p{x} + \lambda f\p{y}\) for each \(x, y \in \calH\) and \(0 \leq \lambda \leq 1\),

\item \emph{\(\mu\)-strongly convex} if \(f\) is proper and \(f-(\mu/2)\norm{\cdot}^{2}\) is convex,

\item \emph{\(\widetilde{\mu}\)-weakly convex} if \(f+(\widetilde{\mu}/2)\norm{\cdot}^{2}\) is convex,

\item \emph{\(L\)-smooth} if \(f\) is Fr\'echet differentiable and the gradient \(\nabla f:\calH \to \calH\) is \(L\)-Lipschitz continuous, i.e., \(\norm{\nabla f(x) - \nabla f(y)} \leq L \norm{x - y}\) for each \(x,y\in\calH\), and

\item \emph{\(\mu_{\textup{gd}}\)-gradient dominated} if \(f\) is Fr\'echet differentiable and
        \begin{align}\label{eq:PL}
            f\p{x} - \inf_{y\in\calH}f\p{y} \leq \frac{1}{2\mu_{\textup{gd}}}\norm{\nabla f(x)}^2
        \end{align}
for each \(x \in \calH\).
Inequality \eqref{eq:PL} is sometimes called the \emph{Polyak--\L{}ojasiewicz} inequality or simply the \emph{\L{}ojasiewicz} inequality.
    \end{enumerate}
\end{definition}

The \emph{Fr\'echet subdifferential} of a function \(f\colon \calH \to \reals \cup \{\pm\infty\}\) is the set-valued operator \(\partial f:\calH \rightrightarrows \calH\) given by
\begin{align*}
    \partial f \p{x} = \begin{cases} \Bigset{u\in\calH \xmiddle| \displaystyle{\liminf_{y\rightarrow x}} \frac{f\p{y}-f\p{x} - \inner{u}{y-x} }{\norm{y-x}}\geq 0 } & \text{if } \abs{f(x)} < + \infty,\\
    \emptyset & \text{otherwise} \end{cases}
\end{align*}
for each \(x\in\calH\).
\begin{enumerate}[(i)]

\item If \(f\) is Fr\'echet differentiable at a point \(x\in\calH\), then \(\partial f \p{x} = \set{\nabla f(x)}\) \cite[Proposition 1.87]{mordukhovich2006variationalanalysisgeneralized}.

\item If \(f\) is proper and convex, the Fr\'echet subdifferential becomes the \emph{convex subdifferential}, i.e., \(\partial f \p{x} = \set{u \in \calH \xmiddle\vert \forall y \in \calH,\, f(y) \geq f(x) + \inner{u}{y-x} }\) for each \(x \in \calH\).

\item If \(f\) is proper and \(\widetilde{\mu}\)-weakly convex for some \(\widetilde{\mu}\in\reals_{++}\), then \(\partial f(x) = \partial \p{f + (\widetilde{\mu}/2)\norm{\cdot}^{2}}\p{x} - \widetilde{\mu} x\) for each \(x\in\mathcal{H}\) \cite[Proposition 1.107 (i)]{mordukhovich2006variationalanalysisgeneralized}, where \(\partial \p{f + (\widetilde{\mu}/2)\norm{\cdot}^{2}}\) reduces to the convex subdifferential.

\end{enumerate}

\begin{definition}\label{def:prox}
Suppose that \(f:\calH\to\reals\cup\set{\pm\infty}\) is proper, lower semicontinuous, and convex, and let \(\gamma\in\reals_{++}\).
Then the \emph{proximal operator} of \(f\) with \emph{step size} \(\gamma\), denoted \(\prox_{\gamma f} : \calH \to \calH \), is defined as the single-valued operator given by
    \begin{align*}
        \prox_{\gamma f}(x) = \argmin_{z\in\calH}\Bigp{f(z) + \frac{1}{2\gamma}\norm{x-z}^2}
    \end{align*}
for each \(x\in\calH\).
\end{definition}

Suppose that \(f:\calH\to\reals\cup\set{\pm\infty}\) is proper, lower semicontinuous, and convex, and let \(\gamma\in\reals_{++}\).
If \(x,p\in\calH\), then \(p = \prox_{\gamma f}(x)\) \(\Leftrightarrow\) \(\gamma^{-1}\p{x-p} \in \partial f (p)\).
Moreover, the \emph{conjugate} of \(f\), denoted \(f^{*}:\calH\to \reals \cup \{+\infty\}\), is the proper, lower semicontinuous and convex function given by \(f^{*}(u) = \sup_{x\in\calH}\p{\inner{u}{x} - f(x)}\) for each \(u\in\calH\).
If \(x,u\in\calH\), then \(u \in \partial f(x) \) \(\Leftrightarrow\) \(x \in \partial f^{*}(u)\).

Let \(G:\calH\rightrightarrows\calH\) be a set-valued operator.
The set of \emph{zeros} of \(G\) is denoted by \(\zer G =\set{x\in\calH \xmiddle| 0 \in G\p{x} }\) and the \emph{graph} of \(G\) is denoted by $\gra G = \set{(x,y)\in\calH\times\calH \xmiddle| y\in G(x)}$.

\begin{definition}\label{def:set_op}
Let \(G:\calH\rightrightarrows\calH\) and \(\mu \in\reals_{++}\).
The operator \(G\) is said to be
    \begin{enumerate}[(i)]
\item \emph{monotone} if \(\inner{u - v}{x-y}\geq 0\) for each \((x,u),(y,v)\in\gra G\),

\item \emph{maximally monotone} if \(G\) is monotone and there does not exist a monotone operator \(H:\calH\rightrightarrows\calH\) such that \(\gra G \subsetneq \gra H\), and

\item \(\mu\)-\emph{strongly monotone} if \(\inner{u - v}{x-y}\geq \mu \norm{x-y}^{2} \) for each \((x,u),(y,v)\in\gra G\).
    \end{enumerate}
\end{definition}



\begin{definition}\label{def:single_op}
Let \(G:\calH \to \calH\), \(L\in\reals_{+}\), and \(\beta\in\reals_{++}\).
The operator \(G\) is said to be
    \begin{enumerate}[(i)]
\item \(L\)-\emph{Lipschitz continuous} if \(\norm{G\p{x} - G\p{y} } \leq L \norm{x-y}\) for each \(x,y\in \calH\), and
\item \(\beta\)-\emph{cocoercive} if \(\inner{G\p{x} - G\p{y} }{x-y} \geq \beta \norm{G\p{x} - G\p{y}}^{2} \) for each \(x,y\in \calH\).
    \end{enumerate}
\end{definition}

We introduce the following conventions that enable us to treat single-valued and singleton-valued operators interchangeably.
\begin{enumerate}[(i)]
\item For notational convenience (at the expense of a slight abuse of notation), we will sometimes identify the operator \(G:\calH\to\calH\) with the set-valued mapping \(\calH \ni x \mapsto \set{G\p{x}} \subseteq\calH \), which will be clear from context.
For example, if \(x,y \in \calH\), the inclusion \(y\in G\p{x}\) should be interpreted as the equality \(y = G\p{x}\).
\item Similarly, if \(G:\calH\rightrightarrows \calH\) and \(T: \calH \to \calH\) satisfy \(G\p{x} = \set{T\p{x}}\) for each \(x\in\calH\), i.e., \(G\) is a singleton-valued operator, we will sometimes identify \(G\) with the corresponding single-valued operator \(T\).
\end{enumerate}

Given any positive integer \(n\), we let the inner-product \(\inner{\cdot}{\cdot}\) on \(\calH^{n}\) be given by
\begin{align*}
    \inner{\bz^{1}}{\bz^{2}}=\sum_{j=1}^{n}\inner{z^{1}_{j}}{z^{2}_{j}}
\end{align*}
for each \(\bz^{i}=\p{z^{i}_{1},\ldots,z^{i}_{n}}\in\calH^{n}\) and \(i\in\llbracket 1,2\rrbracket\).
If \(M\in \reals^{m\times n}\), we define the tensor product \(M\kron \Id\) to be the mapping \((M\kron \Id):\calH^{n}\rightarrow\calH^{m}\) such that
\begin{align*}
  (M\kron\Id)\bz = \Bigp{\sum_{j=1}^n[M]_{1,j}z_{j},\ldots,\sum_{j=1}^n[M]_{m,j}z_{j}}
\end{align*}
for each \(\bz=\p{z_{1},\ldots,z_{n}}\in\calH^n\). 
We define the mapping \(\mathcal{Q}:\sym^{n}\times \calH^{n} \rightarrow \reals\) by \(\quadform{M}{\bz}=\inner{\bz}{(M\kron\Id)\bz}\) for each \(M\in\sym^{n}\) and \(\bz\in\calH^{n}\).

\section{Problem class and algorithm representation}
\label{sec:problem_and_algorithm}

In this section, we introduce the class of optimization and inclusion problems we consider and the representations of algorithms that solve them.

\subsection{Problem class}\label{sec:problem_class}

To cover both structured optimization and inclusion problems, we introduce two disjoint index sets \(\IndexFunc, \IndexOp \subseteq \llbracket1,m\rrbracket\), where \(m\in\mathbb{N}\), such that \(\IndexFunc\cup \IndexOp = \llbracket1,m\rrbracket\), and consider inclusion problems of the form
\begin{align}
\label{eq:the_problem_inclusion} \text{find}\ y\in\calH\ \text{ such that }\ 0\in\sum_{i\in\IndexFunc} \partial f_{i}\p{y} + \sum_{i\in\IndexOp} G_{i}\p{y},
\end{align}
where the functions \(f_i:\calH \to \reals \cup \set{\pm \infty}\) and operators \(G_i : \calH \rightrightarrows \calH\) are chosen from some user-specified function class \(\mathcal{F}_i\) and operator class \(\mathcal{G}_{i}\), respectively, i.e.,
\begin{align*}
    \p{\forall i \in \IndexFunc}& \quad f_{i} \in \mathcal{F}_{i} \subseteq \set{ f: \calH \to \reals\cup \set{\pm \infty}},\\
    \p{\forall i \in \IndexOp}& \quad G_{i} \in \mathcal{G}_{i} \subseteq \set{ G: \calH \rightrightarrows \calH }.
\end{align*}
For example, if \(\IndexOp=\emptyset\), then \eqref{eq:the_problem_inclusion} is a first-order optimality condition for minimizing \(\sum_{i\in\IndexFunc} f_{i}\).
Moreover, \eqref{eq:the_problem_inclusion} provides a formalism that covers monotone inclusion problems \cite{bauschke2017convexanalysismonotone,ryu2022largescaleconvex}, certain equilibrium problems \cite{bricenoarias2013monotoneoperatormethods,bui2021analysisnumericalsolution,combettes2005equilibrium}, so-called (mixed) variational inequalities \cite{facchinei2004finitedimensionalvariationali,facchinei2004finitedimensionalvariationalii,korpelevich1976extragradientmethodfinding}, and beyond.

\subsubsection{Interface and examples}

The \texttt{InclusionProblem} class in the \texttt{autolyap.problemclass} module of \texttt{AutoLyap} provides the interface for formulating inclusion problems of the form \eqref{eq:the_problem_inclusion}.
A few function and operator classes are listed in \Cref{tab:function_classes,tab:operator_classes}, respectively.

For instance, the minimization problem
\begin{align*}
    \minimize_{y\in\calH}\; f_{1}\p{y} + f_{2}\p{y},
\end{align*}
where \(f_1,f_2:\calH\to\reals\cup\set{\pm\infty}\) are proper, lower semicontinuous, and convex, is modeled by the inclusion problem
\begin{align}\label{eq:inclusion_problem:convex_convex}
    \text{find}\ y\in\calH\ \text{ such that }\ 0\in \partial f_{1}\p{y} + \partial f_{2}\p{y},
\end{align}

which in turn is modeled by
\begin{lstlisting}[label={lst:convex_convex}]
from autolyap.problemclass import Convex, InclusionProblem

components_list = [ 
    Convex(), # $\color{gray}{f_1}$
    Convex() # $\color{gray}{f_2}$
] 
problem = InclusionProblem(components_list) # (*@\textcolor{gray}{\eqref{eq:inclusion_problem:convex_convex}}@*)
\end{lstlisting}
in \texttt{AutoLyap}. 

The API also supports intersections of function classes and intersections of operator classes (but not mixed function-operator intersections). 
For example, suppose instead that $\partial f_1$ in \eqref{eq:inclusion_problem:convex_convex} is replaced with an operator \(G_1:\calH\rightarrow\calH\) that is \(\mu\)-strongly monotone and \(L\)-Lipschitz continuous for \(\p{\mu,L} = \p{1,2}\), i.e.,
\begin{align}\label{eq:inclusion_problem:maximally_monotone_plus_Lipschitz_convex}
    \text{find}\ y\in\calH\ \text{ such that }\ 0\in G_{1}\p{y} + \partial f_{2}\p{y}.
\end{align}
This is modeled by
\begin{lstlisting}[label={lst:maximally_monotone_plus_Lipschitz_convex}]
from autolyap.problemclass import (
    InclusionProblem,
    StronglyMonotone,
    LipschitzOperator,
    Convex
)

mu, L = 1, 2

components_list = [ 
    [StronglyMonotone(mu), LipschitzOperator(L)], # $\color{gray}{G_1}$
    Convex() # $\color{gray}{f_2}$
] 
problem = InclusionProblem(components_list) # (*@\textcolor{gray}{\eqref{eq:inclusion_problem:maximally_monotone_plus_Lipschitz_convex}}@*)
\end{lstlisting}
in \texttt{AutoLyap}.

\begin{table}[t]
\centering 
\caption{
Some function classes included in the \texttt{autolyap.problemclass} module of \texttt{AutoLyap}.
See \Cref{def:func_defs} for formal definitions.
Further details are found in the \href{https://autolyap.github.io/function_classes/}{documentation}. 
}
  \label{tab:function_classes}
  \begin{tabular}{l p{0.5\textwidth}}
    \toprule
    \textbf{Class} & \textbf{Description} \\
    \midrule
    \rowcolor{black!10}[\tabcolsep][\tabcolsep]
    \texttt{Convex}
      & Class of functions \(f:\calH\to\reals\cup\set{\pm\infty}\) that are proper, lower semicontinuous, and
      convex.
      \\
    \texttt{StronglyConvex(\(\mu\))}
      & Class of functions \(f:\calH\to\reals\cup\set{\pm\infty}\) that are proper, lower semicontinuous, and
      \(\mu\)-strongly convex for some
      \(\mu\in\reals_{++}\).
      \\
    \rowcolor{black!10}[\tabcolsep][\tabcolsep]
    \texttt{WeaklyConvex(\(\tilde{\mu}\))}
      & Class of functions \(f:\calH\to\reals\cup\set{\pm\infty}\) that are proper, lower semicontinuous, and
      \(\tilde{\mu}\)-weakly convex for some \(\tilde{\mu}\in\reals_{++}\).
      \\
    \texttt{Smooth(\(L\))}
      & Class of functions \(f:\calH\to\reals\) that are \(L\)-smooth for some
      \(L\in\reals_{++}\).
      \\
    \rowcolor{black!10}[\tabcolsep][\tabcolsep]
    \texttt{SmoothConvex(\(L\))}
      & Class of functions \(f:\calH\to\reals\) that are convex and \(L\)-smooth for some
      \(L\in\reals_{++}\).
      \\
    \texttt{SmoothStronglyConvex(\(\mu,L\))}
      & Class of functions \(f:\calH\to\reals\) that are \(\mu\)-strongly convex and
      \(L\)-smooth for some \(\mu,L\in\reals_{++}\) such that \(\mu < L\).
      \\
    \rowcolor{black!10}[\tabcolsep][\tabcolsep]
    \texttt{SmoothWeaklyConvex(\(\tilde{\mu},L\))}
      & Class of functions \(f:\calH\to\reals\) that are \(\tilde{\mu}\)-weakly convex and
      \(L\)-smooth for some \(\tilde{\mu},L\in\reals_{++}\).
      \\
    \texttt{GradientDominated(\(\mu_{\textup{gd}}\))}
      & Class of functions \(f:\calH\to\reals\) that are
      \(\mu_{\textup{gd}}\)-gradient dominated for some
      \(\mu_{\textup{gd}}\in\reals_{++}\). Requires that \(m=1\) and \(\IndexOp = \emptyset\) in
      \eqref{eq:the_problem_inclusion}.
      \\
    \bottomrule
  \end{tabular}
\end{table}

\begin{table}[t]
\centering 
\caption{
Some operator classes included in the \texttt{autolyap.problemclass} module of \texttt{AutoLyap}.
See \Cref{def:set_op,def:single_op} for formal definitions.
Further details are found in the \href{https://autolyap.github.io/operator_classes/}{documentation}.
}
  \label{tab:operator_classes}
  \begin{tabular}{l p{0.6\textwidth}}
    \toprule
    \textbf{Class} & \textbf{Description} \\
    \midrule
    \rowcolor{black!10}[\tabcolsep][\tabcolsep]
    \texttt{MaximallyMonotone}
      & Class of operators \(G:\mathcal H\rightrightarrows\mathcal H\) that are maximally monotone. \\
    \texttt{StronglyMonotone(\(\mu\))}
      & Class of operators \(G:\mathcal H\rightrightarrows\mathcal H\) that are \(\mu\)-strongly and maximally
      monotone for some \(\mu\in\reals_{++}\). \\
    \rowcolor{black!10}[\tabcolsep][\tabcolsep]
    \texttt{LipschitzOperator(\(L\))}
      & Class of operators \(G:\mathcal H\rightarrow\mathcal H\) that are \(L\)-Lipschitz continuous for some
      \(L\in\reals_{++}\). \\
    \texttt{Cocoercive(\(\beta\))}
      & Class of operators \(G:\mathcal H\rightarrow\mathcal H\) that are \(\beta\)-cocoercive for some
      \(\beta\in\reals_{++}\). \\
    \bottomrule
  \end{tabular}
\end{table}

\subsection{Algorithm representation}\label{subsec:state_space}

We consider first-order methods for solving \eqref{eq:the_problem_inclusion} that can be represented as discrete-time linear time-varying state-space systems interconnected in feedback with the potentially nonlinear and set-valued operators \((\partial f_{i})_{i\in\IndexFunc}\) and \((G_{i})_{i\in\IndexOp}\) that define the problem.

Before presenting the algorithm representation, we introduce the necessary notation:
\begin{enumerate}[(i)]
\item \(\NumFunc = \abs{\IndexFunc}\) denotes the number of functional components in \eqref{eq:the_problem_inclusion};
\item \(\NumOp = \abs{\IndexOp}\) denotes the number of operator components in \eqref{eq:the_problem_inclusion};
\item \(\NumEval_{i}\in\mathbb{N}\) denotes the number of evaluations of \(\partial f_{i}\) and \(G_{i}\) per iteration, for \(i\in\IndexFunc\) and \(i\in\IndexOp\), respectively;
\item \(\NumEvalFunc = \sum_{i\in\IndexFunc} \NumEval_{i}\) denotes the total number of subdifferential evaluations per iteration;
\item \(\NumEvalOp = \sum_{i\in\IndexOp} \NumEval_{i}\) denotes the total number of operator evaluations per iteration; and
\item \(\NumEval = \NumEvalFunc + \NumEvalOp\) denotes the combined total number of evaluations per iteration.
\end{enumerate}

Since we consider methods that allow for multiple evaluations of \((\partial f_{i})_{i\in\IndexFunc}\) and \((G_{i})_{i\in\IndexOp}\) per iteration, we define \(\bfcn_{i}:\calH^{\NumEval_{i}}\to\p{\reals\cup\{\pm\infty\}}^{\NumEval_{i}}\) such that
\begin{align*}
    \Bigp{ \begin{array}{@{}c@{}} \forall i\in\IndexFunc\\
    \forall \by_{i}=\p{y_{i,1},\ldots,y_{i,\NumEval_{i}}}\in \calH^{\NumEval_{i}} \end{array} } \quad \bfcn_{i}(\by_{i})&=\p{f_i\p{y_{i,1}},\ldots,f_{i}\p{y_{i,\NumEval_{i}}}},
\end{align*}
\(\bm{\partial}\bfcn_{i}:\calH^{\NumEval_{i}}\rightrightarrows\calH^{\NumEval_{i}}\) such that\footnote{In this context, the symbol \(\Pi\) is used for Cartesian products.}
\begin{align*}
    \Bigp{ \begin{array}{@{}c@{}} \forall i\in\IndexFunc\\
    \forall \by_{i}=\p{y_{i,1},\ldots,y_{i,\NumEval_{i}}}\in \calH^{\NumEval_{i}} \end{array} } \quad \bm{\partial}\bfcn_{i}(\by_{i})&= \prod_{j=1}^{\NumEval_{i}} \partial f_{i}\p{y_{i,j}},
\end{align*}
and \(\bm{G}_{i}: \calH^{\NumEval_{i}} \rightrightarrows \calH^{\NumEval_{i}}\) such that
\begin{align*}
    \Bigp{ \begin{array}{@{}c@{}} \forall i\in\IndexOp\\
    \forall \by_{i}=\p{y_{i,1},\ldots,y_{i,\NumEval_{i}}}\in \calH^{\NumEval_{i}} \end{array} } \quad \bm{G}_{i}\p{\by_{i}}=\prod_{j=1}^{\NumEval_{i}} G_{i}\p{y_{i,j}}.
\end{align*}

We are now ready to give the algorithm representation: Pick an initial \(\bx_{0}\in\calH^{n}\), an iteration horizon \(K\in\naturals\cup\set{\infty}\), and let
\begin{equation}\label{eq:linear_system_with_nonlinearity}
    \p{\forall k \in\llbracket0,K\rrbracket} \quad \left[ \begin{aligned} &\bx^{k+1} = \p{A_{k} \kron \Id} \bx^{k} + \p{B_{k} \kron \Id} \bu^{k},\\
    &\by^{k} = \p{C_{k} \kron \Id} \bx^{k} + \p{D_{k} \kron \Id} \bu^{k},\\
    &\p{\bu^{k}_{i}}_{i\in\IndexFunc} \in \prod_{i\in\IndexFunc}\bm{\partial}\bfcn_{i}\p{\by_{i}^{k}},\\
    &\p{\bu^{k}_{i}}_{i\in\IndexOp} \in \prod_{i\in\IndexOp}\bm{G}_{i}\p{\by_{i}^{k}},\\
    &\bFcn^{k} =\p{\bfcn_{i}\p{\by_{i}^{k}} }_{i\in\IndexFunc}, \end{aligned} \right.
\end{equation}
where
\begin{align*}
    \bx^{k}&=\p{x^{k}_{1},\ldots,x^{k}_{n}}\in\calH^n, \\
    \bu^{k}&=\p{\bu^{k}_{1},\ldots,\bu^{k}_{m}}\in\prod_{i=1}^{m}\calH^{\NumEval_{i}}, \quad \bu^{k}_{i} = \p{u^{k}_{i,1},\ldots,u^{k}_{i,\NumEval_{i}}} \in \calH^{\NumEval_{i}},\\    
    \by^{k}&=\p{\by^{k}_{1},\ldots,\by^{k}_{m}}\in\prod_{i=1}^{m}\calH^{\NumEval_{i}}, \quad \by^{k}_{i} = \p{y^{k}_{i,1},\ldots,y^{k}_{i,\NumEval_{i}}} \in \calH^{\NumEval_{i}},\\
    \bFcn^{k}&\in\mathbb{R}^{\NumEvalFunc},
\end{align*}
are the algorithm variables and
\begin{align}\label{eq:ABCD}
    A_{k}&\in\reals^{n\times n},& B_{k}&\in\reals^{n\times {\NumEval}},& C_{k}&\in\reals^{{\NumEval}\times n},& D_{k}&\in\reals^{{\NumEval}\times {\NumEval}}
\end{align}
are matrices containing the parameters of the method at hand.

\subsubsection{Interface and examples}

An interface for \eqref{eq:linear_system_with_nonlinearity} is provided by the abstract base class \texttt{Algorithm} and its concrete subclasses, which can be found in the \texttt{autolyap.algorithms} module of \texttt{AutoLyap}. Specifically, each algorithm in \texttt{AutoLyap} must be implemented as a concrete subclass of \texttt{Algorithm}, and must define the abstract method \texttt{get\_ABCD}. Some concrete classes included in \texttt{AutoLyap} are summarized in \Cref{tab:concrete_algorithms}.

\begin{table}[t]
\centering
\caption{
Some concrete algorithm classes included in the \texttt{autolyap.algorithms} module of \texttt{AutoLyap}.
Further details are found in the \href{https://autolyap.github.io/concrete_algorithms/}{documentation}.
}
  \label{tab:concrete_algorithms}
  \begin{tabular}{l p{0.62\textwidth}}
    \toprule
    \textbf{Class} & \textbf{Description} \\
    \midrule
    \rowcolor{black!10}[\tabcolsep][\tabcolsep]
    \texttt{AcceleratedProximalPoint}
      & An improved proximal point method~\cite{kim2021acceleratedproximalpoint}. \\
    \texttt{ChambollePock}
      & Chambolle--Pock primal-dual method \cite{chambolle2011firstorderprimal}. \\
    \rowcolor{black!10}[\tabcolsep][\tabcolsep]
    \texttt{DavisYin}
      & Davis--Yin three-operator splitting \cite{davis2017threeoperatorsplitting}. \\
    \texttt{DouglasRachford}
      & Douglas--Rachford splitting \cite{douglas1956numericalsolutionheat,eckstein1992douglasrachfordsplittingmethod,lions1979splittingalgorithmssum}. \\
    \rowcolor{black!10}[\tabcolsep][\tabcolsep]
    \texttt{Extragradient}
      & Extragradient method \cite{korpelevich1976extragradientmethodfinding}. \\
    \texttt{ForwardMethod}
      & Forward method. \\
    \rowcolor{black!10}[\tabcolsep][\tabcolsep]
    \texttt{GradientMethod}
      & Gradient method. \\
    \texttt{GradientNesterovMomentum}
      & Gradient method with Nesterov-like momentum \cite{nesterov2018lecturesconvexoptimization}. \\
    \rowcolor{black!10}[\tabcolsep][\tabcolsep]
    \texttt{HeavyBallMethod}
      & Heavy-ball method \cite{Polyak1964HeavyBall}. \\
    \texttt{ITEM}
      & Information-theoretic exact method (ITEM) \cite{taylor2023optimalgradientmethod}. \\
    \rowcolor{black!10}[\tabcolsep][\tabcolsep]
    \texttt{MalitskyTamFRB}
      & Malitsky--Tam forward-reflected-backward method \cite{malitsky2020forwardbackwardsplitting}. \\
    \texttt{NesterovConstant}
      & Nesterov's constant-step scheme \cite{nesterov2018lecturesconvexoptimization}. \\
    \rowcolor{black!10}[\tabcolsep][\tabcolsep]
    \texttt{NesterovFastGradientMethod}
      & Nesterov's fast gradient method \cite{nesterov1983fast}. \\
    \texttt{OptimizedGradientMethod}
      & Optimized gradient method \cite{kim2015optimizedfirstorder}. \\
    \rowcolor{black!10}[\tabcolsep][\tabcolsep]
    \texttt{ProximalPoint}
      & Proximal point method \cite{Moreau1965,Martinet1970,Rockafellar1976PPA}. \\
    \texttt{TripleMomentum}
      & Triple-momentum method \cite{vanscoy2018fastestknownglobally}. \\
    \rowcolor{black!10}[\tabcolsep][\tabcolsep]
    \texttt{TsengFBF}
      & Tseng's forward-backward-forward method \cite{tseng2000modifiedforwardbackward}. \\
    \bottomrule
  \end{tabular}
\end{table}

A simplified version of the \texttt{Algorithm} interface is shown below:
\begin{lstlisting}[label={lst:class_Alg_def}]
from abc import ABC, abstractmethod
from typing import List, Tuple
import numpy as np

class Algorithm(ABC):
    def __init__(
        self,
        n: int,                    # dimension $\color{gray}{n}$ of $\color{gray}{\bx^{k}}$ in (*@\textcolor{gray}{\eqref{eq:linear_system_with_nonlinearity}}@*)
        m: int,                    # number of components $\color{gray}{m}$ in (*@\textcolor{gray}{\eqref{eq:the_problem_inclusion}}@*)
        m_bar_is: List[int],       # list of $\color{gray}{(\NumEval_{i})_{i=1}^m}$ in (*@\textcolor{gray}{\eqref{eq:linear_system_with_nonlinearity}}@*)
        I_func: List[int],         # index set $\color{gray}{\IndexFunc}$ in (*@\textcolor{gray}{\eqref{eq:the_problem_inclusion}}@*)
        I_op: List[int],           # index set $\color{gray}{\IndexOp}$ in (*@\textcolor{gray}{\eqref{eq:the_problem_inclusion}}@*)
    ) -> None:
        # ... (omitted for brevity) ...

    @abstractmethod
    def get_ABCD(
        self,
        k: int,                    # iteration $\color{gray}{k}$ in (*@\textcolor{gray}{\eqref{eq:linear_system_with_nonlinearity}}@*)
    ) -> Tuple[
        np.ndarray,                # $\color{gray}{A_k}$ in (*@\textcolor{gray}{\eqref{eq:ABCD}}@*)
        np.ndarray,                # $\color{gray}{B_k}$ in (*@\textcolor{gray}{\eqref{eq:ABCD}}@*)
        np.ndarray,                # $\color{gray}{C_k}$ in (*@\textcolor{gray}{\eqref{eq:ABCD}}@*)
        np.ndarray,                # $\color{gray}{D_k}$ in (*@\textcolor{gray}{\eqref{eq:ABCD}}@*)
    ]:
        pass
\end{lstlisting}

Below, we provide two illustrative examples showing how to convert an algorithm from standard form into the algorithm representation \eqref{eq:linear_system_with_nonlinearity}, and then implement it as a concrete subclass of \texttt{Algorithm}.
Since the algorithm representation is generally not unique, we typically aim to choose one with the smallest possible state dimension \(n\).

\paragraph{The Chambolle--Pock method. }

Starting from \eqref{eq:inclusion_problem:convex_convex}, suppose that \(f_1,f_2:\calH\to\reals\cup\set{\pm\infty}\) are proper, lower semicontinuous, and convex.
In the identity-operator case, the Chambolle--Pock method \cite{chambolle2011firstorderprimal}, also known as the primal-dual hybrid gradient (PDHG) method, is given by
\begin{equation}\label{eq:chambolle_pock_interface:standard}
    \p{\forall k \in\naturals} \quad \left[ \begin{aligned}
    x^{k+1} &= \prox_{\tau f_1}\p{x^{k}-\tau y^{k}},\\
    y^{k+1} &= \prox_{\sigma f_2^{*}}\p{y^{k}+\sigma\p{x^{k+1} + \theta \p{x^{k+1} - x^{k}} }},
    \end{aligned} \right.
\end{equation}
where \(\tau,\sigma\in\reals_{++}\) are primal and dual step sizes, respectively, \(\theta\in\reals\) is a relaxation parameter, and \(f_2^*\) is the conjugate of \(f_2\).

Using the proximal-operator optimality condition, we obtain
\begin{equation*}
    \begin{aligned}
    \frac{x^{k}-\tau y^{k} - x^{k+1}}{\tau} &\in \partial f_1\p{x^{k+1}},\\
    \frac{y^{k}+\sigma\p{x^{k+1} + \theta \p{x^{k+1} - x^{k}} } - y^{k+1}}{\sigma} &\in \partial f_2^{*}\p{y^{k+1}},
    \end{aligned}
\end{equation*}
or equivalently, 
\begin{equation*}
    \begin{aligned}
    \frac{x^{k}-\tau y^{k} - x^{k+1}}{\tau} &\in \partial f_1\p{x^{k+1}},\\
    y^{k+1} &\in \partial f_2\left(\frac{y^{k}+\sigma\p{x^{k+1} + \theta \p{x^{k+1} - x^{k}} } - y^{k+1}}{\sigma}\right),
    \end{aligned}
\end{equation*}
which can be written as
\begin{equation*}
    \begin{aligned}
    x^{k+1} &\in x^{k} - \tau y^{k}  - \tau \partial f_1\p{x^{k+1}},\\
    y^{k+1} &\in \partial f_2\left(\frac{y^{k}+\sigma\p{x^{k+1} + \theta \p{x^{k+1} - x^{k}} } - y^{k+1}}{\sigma}\right).
    \end{aligned}
\end{equation*}

Picking
\begin{align*}
    \bx^{k} &= \p{x^{k},y^{k}},\\
    \bu^{k} &= \left(\frac{x^{k}-\tau y^{k}-x^{k+1}}{\tau}, y^{k+1}\right),\\
    \by^{k} &= \left(x^{k+1}, \frac{y^{k}+\sigma\p{x^{k+1} + \theta \p{x^{k+1} - x^{k}} } - y^{k+1}}{\sigma}\right), \\
    \bm{\partial}\bfcn\p{\by} &= \partial f_1\p{y_1}\times\partial f_2\p{y_2} \text{ for each } \by=\p{y_1,y_2}\in\calH^2,
\end{align*}
we see that
\begin{equation}\label{eq:chambolle_pock_interface:state_space}
    \begin{aligned}
    \bx^{k+1} &= \Bigp{\begin{bmatrix}
        1 & -\tau\\
        0 & 0
    \end{bmatrix}\kron\Id}\bx^k + \Bigp{\begin{bmatrix}
        -\tau & 0\\
        0 & 1
    \end{bmatrix}\kron\Id}\bu^k,\\
    \by^{k} &= \Bigp{\begin{bmatrix}
        1 & -\tau\\
        1 & \frac{1}{\sigma} - \tau\p{1+\theta}
    \end{bmatrix}\kron\Id}\bx^k + \Bigp{\begin{bmatrix}
        -\tau & 0\\
        -\tau\p{1+\theta} & -\frac{1}{\sigma}
    \end{bmatrix}\kron\Id}\bu^k,\\
    \bu^{k} &\in \bm{\partial}\bfcn\p{\by^{k}},
    \end{aligned}
\end{equation}
which fits algorithm representation \eqref{eq:linear_system_with_nonlinearity}. Moreover, the structural parameters are
\begin{align*}
    n = 2,\quad m = 2,\quad \p{\NumEval_{i}}_{i=1}^{m} = \p{1,1},\quad
    \IndexFunc = \set{1,2},\quad \IndexOp = \emptyset,
\end{align*}
and in \texttt{AutoLyap}, this is implemented as:

\begin{lstlisting}[label={lst:chambolle_pock_impl}]
class ChambollePock(Algorithm):
    def __init__(self, tau: float, sigma: float, theta: float) -> None:
        super().__init__(2, 2, [1, 1], [1, 2], [])
        self.tau = tau
        self.sigma = sigma
        self.theta = theta

    def get_ABCD(self, k: int) -> Tuple[np.ndarray, np.ndarray, np.ndarray, np.ndarray]:
        A = np.array([[1, -self.tau],
                      [0, 0]])
        B = np.array([[-self.tau, 0],
                      [0, 1]])
        C = np.array([[1, -self.tau],
                      [1, 1 / self.sigma - self.tau * (1 + self.theta)]])
        D = np.array([[-self.tau, 0],
                      [-self.tau * (1 + self.theta), -1 / self.sigma]])
        return (A, B, C, D)
\end{lstlisting}

This is exactly how \texttt{ChambollePock} is implemented in \texttt{AutoLyap}. It can be imported with:
\begin{lstlisting}[label={lst:chambolle_pock_import}]
from autolyap.algorithms import ChambollePock
\end{lstlisting}

\paragraph{Nesterov's fast gradient method. }

Suppose that \(f_{1}:\calH\to\reals\) is convex and \(L\)-smooth for some \(L\in\reals_{++}\), and let \(\gamma\in\reals_{++}\), \(\lambda_{0}=1\), and \(x^{-1},x^{0}\in\calH\).
Nesterov's fast gradient method \cite{nesterov1983fast} is given by
\begin{align}\label{eq:nesterov_fast_interface:standard_form}
    \p{\forall k \in \naturals} & \quad \left[ \begin{aligned} &y^{k} = x^{k} + \delta_{k} \p{x^{k} - x^{k-1}},\\
    &x^{k+1} = y^{k} - \gamma \nabla f_{1}\p{y^{k}}, \end{aligned} \right.
\end{align}
where
\begin{align*}
    \p{\forall k \in \naturals} \quad \left[ \begin{aligned} &\delta_{k} = \frac{\lambda_{k}-1}{\lambda_{k+1}},\\
    &\lambda_{k+1} = \frac{1 + \sqrt{1 + 4\lambda_{k}^{2} }}{2}. \end{aligned} \right.
\end{align*}
In \texttt{AutoLyap}, we include an additional evaluation \(\nabla f_{1}\p{x^k}\), which is useful for Lyapunov templates involving function values at \(x^k\).

Picking
\begin{align*}
    \bx^{k} &= \p{x^{k},x^{k-1}},\\
    \bu^{k} &= \p{\nabla f_{1}\p{y^{k}}, \nabla f_{1}\p{x^{k}}},\\
    \by^{k} &= \p{y^{k},x^{k}},\\
    \bm{\partial}\bfcn\p{\by} &= \set{\nabla f_{1}\p{y_1}}\times\set{\nabla f_{1}\p{y_2}} \text{ for each } \by=\p{y_1,y_2}\in\calH^2,
\end{align*}
the update \eqref{eq:nesterov_fast_interface:standard_form} can be written as
\begin{equation}\label{eq:nesterov_fast_interface:state_space}
    \begin{aligned}
    \bx^{k+1} &= \Bigp{\begin{bmatrix}
        1+\delta_k & -\delta_k\\
        1 & 0
    \end{bmatrix}\kron\Id}\bx^k + \Bigp{\begin{bmatrix}
        -\gamma & 0\\
        0 & 0
    \end{bmatrix}\kron\Id}\bu^k,\\
    \by^{k} &= \Bigp{\begin{bmatrix}
        1+\delta_k & -\delta_k\\
        1 & 0
    \end{bmatrix}\kron\Id}\bx^k + \Bigp{\begin{bmatrix}
        0 & 0\\
        0 & 0
    \end{bmatrix}\kron\Id}\bu^k,\\
    \bu^{k} &\in \bm{\partial}\bfcn\p{\by^{k}},
    \end{aligned}
\end{equation}
which fits algorithm representation \eqref{eq:linear_system_with_nonlinearity}. Moreover, the structural parameters are
\begin{align*}
    n = 2,\quad m = 1,\quad \p{\NumEval_{i}}_{i=1}^{m} = \p{2},\quad
    \IndexFunc = \set{1},\quad \IndexOp = \emptyset,
\end{align*}
and in \texttt{AutoLyap}, this is implemented as:

\begin{lstlisting}[label={lst:nesterov_fast_impl}]
class NesterovFastGradientMethod(Algorithm):
    def __init__(self, gamma: float) -> None:
        super().__init__(2, 1, [2], [1], [])
        self.gamma = gamma

    def get_ABCD(self, k: int) -> Tuple[np.ndarray, np.ndarray, np.ndarray, np.ndarray]:
        lambda_var = 1
        for _ in range(0, k + 1):
            lambda_var_prev = lambda_var
            lambda_var = (1 + np.sqrt(1 + 4 * lambda_var ** 2)) / 2
        delta = (lambda_var_prev - 1) / lambda_var

        A = np.array([[1 + delta, -delta],
                      [1, 0]])
        B = np.array([[-self.gamma, 0],
                      [0, 0]])
        C = np.array([[1 + delta, -delta],
                      [1, 0]])
        D = np.array([[0, 0],
                      [0, 0]])
        return (A, B, C, D)
\end{lstlisting}

This is exactly how \texttt{NesterovFastGradientMethod} is implemented in \texttt{AutoLyap}. It can be imported with:
\begin{lstlisting}[label={lst:nesterov_fast_import}]
from autolyap.algorithms import NesterovFastGradientMethod
\end{lstlisting}

\section{The analysis tools}
\label{sec:analysis}

In this section, we introduce two classes of Lyapunov analyses.
The first class comprises stationary, or iteration-independent, analyses, while the second involves nonstationary, or iteration-dependent, analyses.
These analyses utilize quadratic ansatzes, facilitating verification through semidefinite programs, as detailed in \Cref{sec:iteration_independent_Lyapunov,sec:iteration_dependent_Lyapunov}, respectively. For later convenience across these sections, we introduce the solution variables used in the Lyapunov analyses. In particular, we are interested in Lyapunov analyses that may depend on a solution to the inclusion problem~\eqref{eq:the_problem_inclusion}. Without loss of generality and for computational efficiency, we will consider the variables
\begin{align}\label{eq:solution}
    \p{y^{\star}, \hat{\bu}^{\star}, \bFcn^{\star}} \in \Bigset{ \underbracket{\p{y,\hat{\bu}, \bFcn}}_{ \in \calH \times\calH^{m-1} \times \reals^{\NumFunc} } \xmiddle| \begin{aligned} &\p{u_{i}}_{i\in\IndexFunc} \in \prod_{i\in\IndexFunc}\partial f_{i}\p{y},\\
    &\p{u_{i}}_{i\in\IndexOp} \in \prod_{i\in\IndexOp}G_{i}\p{y},\\
    &\sum_{i=1}^{m} u_{i} = 0,\\
    & \hat{\bu} = \p{u_{1},\ldots,u_{m-1}},\\
    & \bFcn=\p{\bfcn_{i}\p{y}}_{i\in\IndexFunc} \end{aligned} },
\end{align}
where \(\hat{\bu}^{\star}\) is void when \(m=1\).
For example, it is clear that \(y^{\star}\) in \eqref{eq:solution} is a solution to the inclusion problem~\eqref{eq:the_problem_inclusion}.

\subsection{Iteration-independent Lyapunov analyses}\label{sec:iteration_independent_Lyapunov}
In this case, we consider algorithms that continue to iterate indefinitely and have iteration-independent parameters.
That is, in~\eqref{eq:linear_system_with_nonlinearity}, we assume that \(K=\infty\) and that there exist fixed matrices
\begin{align*}
    \p{A,B,C,D}\in\reals^{n\times n}\times\reals^{n\times \NumEval}\times\reals^{\NumEval\times n}\times\reals^{\NumEval \times \NumEval}
\end{align*}
such that
\begin{align}\label{eq:constant_ABCD}
    \p{\forall k \in \naturals }\quad \p{A_{k},B_{k},C_{k},D_{k}} = \p{A,B,C,D}.
\end{align}

\begin{definition}\label{def:iteration_independent_Lyapunov}
Suppose that \eqref{eq:constant_ABCD} holds, \(\bx_{0}\in\calH^{n}\) is an initial point, \(( ( \bx^{k},\allowbreak \bu^{k},\allowbreak \by^{k},\allowbreak \bFcn^{k} ) )_{k\in\naturals}\) is a sequence of iterates satisfying~\eqref{eq:linear_system_with_nonlinearity}, \(\p{y^{\star},\allowbreak \hat{\bu}^{\star},\allowbreak \bFcn^{\star}}\) is a point satisfying~\eqref{eq:solution}, \(\rho \in [0,1]\) is a contraction factor, \(h\in\naturals\) is a history parameter, and \(\alpha\in\naturals\) is an overlap parameter.
Define
    \begin{equation}
    \begin{aligned}\label{eq:iteration_independent_Lyapunov:V} \p{\forall k \in \naturals} \quad \mathcal{V}(W,w,k) {}={}&{} \quadform{W}{\p{\bx^{k},\bu^{k},\ldots,\bu^{k+h},\hat{\bu}^{\star},y^{\star}}}\\
    & + w^{\top}\p{\bFcn^{k},\ldots,\bFcn^{k+h},\bFcn^{\star}}, \end{aligned}
    \end{equation}
for each \(\p{W,w}\in\set{\p{Q,q},\p{P,p}}\), where
\begin{align*}
    Q,P \in \sym^{n+\p{h+1}\NumEval+m}
    \quad \text{ and } \quad
    q,p\in\mathbb{R}^{\p{h+1}\NumEvalFunc + \NumFunc},
\end{align*}
and define
    \begin{equation}
    \begin{aligned}\label{eq:iteration_independent_Lyapunov:R} \p{\forall k \in \naturals} \quad \mathcal{R}(W,w,k) {}={}&{} \quadform{W}{\p{\bx^{k},\bu^{k},\ldots,\bu^{k+h + \alpha + 1},\hat{\bu}^{\star},y^{\star}}}\\
    &{} + w^{\top}\p{\bFcn^{k},\ldots,\bFcn^{k+h + \alpha + 1},\bFcn^{\star}}, \end{aligned}
    \end{equation}

for each \(\p{W,w}\in\set{\p{S,s},\p{T,t}}\), where 
\begin{align*}
    S,T\in\sym^{n+\p{h+\alpha+2}\NumEval+m}
    \quad \text{ and } \quad
    s, \allowbreak t \allowbreak \in\mathbb{R}^{\p{h+\alpha+2}\NumEvalFunc + \NumFunc}.
\end{align*}
We say that \(\p{Q,q,S,s}\) satisfies a
\begin{align*}
    \p{P,p,T,t,\rho, h, \alpha}\textup{-quadratic Lyapunov inequality}
\end{align*}
for algorithm \eqref{eq:linear_system_with_nonlinearity} over the problem class defined by \((\mathcal{F}_i)_{i\in\IndexFunc}\) and \((\mathcal{G}_i)_{i\in\IndexOp}\) if
        \begin{align}
            \p{\forall k \in \naturals}& \quad \mathcal{V}\p{Q,q,k+\alpha+1} \leq \rho \mathcal{V}\p{Q,q,k} - \mathcal{R}\p{S,s,k}, \tag{C1} \label{eq:C1}\\
            \p{\forall k \in \naturals}& \quad \mathcal{V}\p{Q,q,k} \geq \mathcal{V}\p{P,p,k} \geq 0, \tag{C2} \label{eq:C2}\\
            \p{\forall k \in \naturals}& \quad \mathcal{R}\p{S,s,k} \geq \mathcal{R}\p{T,t,k}\geq 0, \tag{C3} \label{eq:C3}\\
            \p{\forall k \in \naturals}& \quad \mathcal{R}\p{S,s,k+1} \leq \mathcal{R}\p{S,s,k}, \tag{C4} \label{eq:C4}
        \end{align}
hold for each initial point \(\bx_{0}\), for each sequence of iterates \(( ( \bx^{k}, \bu^{k}, \by^{k}, \bFcn^{k} ) )_{k\in\naturals}\), for each point \(\p{y^{\star}, \hat{\bu}^{\star},\bFcn^{\star}}\), for each \(\p{f_{i}}_{i\in\IndexFunc} \in \prod_{i\in\IndexFunc} \mathcal{F}_{i}\), and for each \(\p{G_i}_{i\in\IndexOp} \in \prod_{i\in\IndexOp} \mathcal{G}_i\), where \eqref{eq:C4} is an optional requirement that may be removed.
\end{definition}

In the proposed methodology, the user specifies $\p{P,p,T,t,\rho,h,\alpha}$ and \texttt{AutoLyap} searches for $\p{Q,q,S,s}$ satisfying \eqref{eq:C1}-\eqref{eq:C3} (and optionally \eqref{eq:C4}).
When such a $\p{Q,q,S,s}$ exists, the choice of $\p{P,p,T,t,\rho,h,\alpha}$ determines which convergence properties are implied by \Cref{def:iteration_independent_Lyapunov}.
\begin{itemize}
\item If \(\rho \in [0,1[\), then
    \begin{align*}
        0 \leq \mathcal{V}\p{P,p,k} \leq \mathcal{V}\p{Q,q,k} \leq \rho^{\lfloor k/\p{\alpha+1}\rfloor }\max_{i\in\llbracket0,\alpha\rrbracket}\mathcal{V}\p{Q,q,i} \xrightarrow[k\rightarrow \infty]{} 0,
    \end{align*}
i.e.,
\begin{align*}
    (\mathcal{V}\p{P,p,k})_{k\in\naturals} \text{ converges }    \sqrt[\alpha + 1]{\rho}\text{-linearly to zero}.
\end{align*}
\item If \(\rho = 1\), then
    \begin{align*}
        \p{\forall k \in \naturals}\quad \sum_{i=0}^{k} \mathcal{R}\p{T,t,i} \leq \sum_{i=0}^{k} \mathcal{R}\p{S,s,i} \leq \sum_{j=0}^{\alpha} \mathcal{V}\p{Q,q,j},
    \end{align*}
using a telescoping summation argument.
In particular,
\begin{align*}
   (\mathcal{R}\p{T,t,k})_{k\in\naturals} \text{ is summable},
\end{align*}
converges to zero, and, e.g., 
\begin{align*}
    \min_{i\in \llbracket 0,k \rrbracket } \mathcal{R}\p{T,t,i} \in \mathcal{O}\p{1/k} \text{ as }  k\to\infty.
\end{align*}
If the optional requirement \eqref{eq:C4} holds, we conclude the stronger last-iterate convergence result
\begin{align*}
    \mathcal{R}\p{T,t,k} \in o\p{1/k} \text{ as }  k\to\infty.
\end{align*}
\end{itemize}

\subsubsection{Interface and examples}

In \texttt{AutoLyap}, the interface for \Cref{def:iteration_independent_Lyapunov} is \texttt{IterationIndependent.search\_lyapunov}.
For fixed \(\p{P,p,T,t,\rho,h,\alpha}\), it formulates the SDP feasibility problem associated with \Cref{def:iteration_independent_Lyapunov} and solves it using the backend selected through \texttt{SolverOptions}.

To support this workflow, \texttt{AutoLyap} provides helper functions in \path{IterationIndependent.LinearConvergence} and \path{IterationIndependent.SublinearConvergence} for choosing \(\p{P,p,T,t}\).
These choices fix the target lower bounds \(\mathcal{V}\p{P,p,k}\) and \(\mathcal{R}\p{T,t,k}\).
\Cref{tab:iter_indep_linear_helpers,tab:iter_indep_sublinear_helpers} summarize these built-in \path{get_parameters_*} helper functions.
For linear-convergence analyses, \path{IterationIndependent.LinearConvergence.bisection_search_rho} can automatically search for the smallest feasible contraction factor \(\rho\) for fixed \(\p{P,p,T,t,h,\alpha}\).

\begin{table}[t]
\centering
\caption{
Some helper functions in \texttt{IterationIndependent.LinearConvergence}, with inputs, allowed values, and induced targets \(\mathcal{V}\p{P,p,k}\) and \(\mathcal{R}\p{T,t,k}\).
Notation follows \eqref{eq:linear_system_with_nonlinearity} and \Cref{def:iteration_independent_Lyapunov}, with \(\mathcal{V}\) and \(\mathcal{R}\) defined in \eqref{eq:iteration_independent_Lyapunov:V} and \eqref{eq:iteration_independent_Lyapunov:R}.
When \(\NumFunc = 0\), the vectors \(p\) and \(t\) are omitted (equivalently, helpers return \(\p{P,T}\) instead of \(\p{P,p,T,t}\)).
Further details are found in the \href{https://autolyap.github.io}{documentation}.
}
\label{tab:iter_indep_linear_helpers}
\begin{tabular}{p{0.93\textwidth}}
\toprule
\rowcolor{black!10}[\tabcolsep][\tabcolsep]
\parbox{\linewidth}{
\emph{Assignment:}
\begin{equation*}
\p{P,p,T,t}=\texttt{get\_parameters\_distance\_to\_solution}(\texttt{algo},h,\alpha,i,j,\tau).
\end{equation*}
\emph{Arguments and allowed ranges:}
\begin{center}
\begin{minipage}{0.96\linewidth}
\centering
\texttt{algo}: \texttt{Algorithm};\quad
\(h,\alpha\in\naturals\);\quad
\(i\in\llbracket 1,m\rrbracket\);\quad
\(j\in\llbracket 1,\NumEval_{i}\rrbracket\);\quad
\(\tau\in\llbracket 0,h\rrbracket\).
\end{minipage}
\end{center}
\emph{Induced targets:}
\begin{equation*}
\mathcal{V}\p{P,p,k}=\norm{y_{i,j}^{k+\tau}-y^{\star}}^{2},
\qquad
\mathcal{R}\p{T,t,k}=0.
\end{equation*}
}
\\
\addlinespace[0.6em]
\parbox{\linewidth}{
\emph{Requires:} \(m=\NumFunc=1\).\\
\emph{Assignment:}
\begin{equation*}
\p{P,p,T,t}=\texttt{get\_parameters\_function\_value\_suboptimality}(\texttt{algo},h,\alpha,j,\tau).
\end{equation*}
\emph{Arguments and allowed ranges:}
\begin{center}
\begin{minipage}{0.96\linewidth}
\centering
\texttt{algo}: \texttt{Algorithm};\quad
\(h,\alpha\in\naturals\);\quad
\(j\in\llbracket 1,\NumEval_{1}\rrbracket\);\quad
\(\tau\in\llbracket 0,h\rrbracket\).
\end{minipage}
\end{center}
\emph{Induced targets:}
\begin{equation*}
\mathcal{V}\p{P,p,k}=f_{1}\p{y_{1,j}^{k+\tau}}-f_{1}\p{y^{\star}},
\qquad
\mathcal{R}\p{T,t,k}=0.
\end{equation*}
}
\\
\bottomrule
\end{tabular}
\end{table}

\begin{table}[t]
\centering
\caption{
Some helper functions in \texttt{IterationIndependent.SublinearConvergence}, with inputs, allowed values, and induced targets \(\mathcal{V}\p{P,p,k}\) and \(\mathcal{R}\p{T,t,k}\).
Notation follows \eqref{eq:linear_system_with_nonlinearity} and \Cref{def:iteration_independent_Lyapunov}, with \(\mathcal{V}\) and \(\mathcal{R}\) defined in \eqref{eq:iteration_independent_Lyapunov:V} and \eqref{eq:iteration_independent_Lyapunov:R}.
When \(\NumFunc = 0\), the vectors \(p\) and \(t\) are omitted (equivalently, helpers return \(\p{P,T}\) instead of \(\p{P,p,T,t}\)).
Further details are found in the \href{https://autolyap.github.io}{documentation}.
}
\label{tab:iter_indep_sublinear_helpers}
\begin{tabular}{p{0.93\textwidth}}
\toprule
\rowcolor{black!10}[\tabcolsep][\tabcolsep]
\parbox{\linewidth}{
\emph{Assignment:}
\begin{equation*}
\p{P,p,T,t}=\texttt{get\_parameters\_fixed\_point\_residual}(\texttt{algo},h,\alpha,\tau).
\end{equation*}
\emph{Arguments and allowed ranges:}
\begin{center}
\begin{minipage}{0.96\linewidth}
\centering
\texttt{algo}: \texttt{Algorithm};\quad
\(h,\alpha\in\naturals\);\quad
\(\tau\in\llbracket 0,h+\alpha+1\rrbracket\).
\end{minipage}
\end{center}
\emph{Induced targets:}
\begin{equation*}
\mathcal{V}\p{P,p,k}=0,
\qquad
\mathcal{R}\p{T,t,k}=\norm{\bx^{k+\tau+1}-\bx^{k+\tau}}^{2}.
\end{equation*}
}
\\
\addlinespace[0.6em]
\parbox{\linewidth}{
\emph{Requires:} \(m=\NumFunc=1\).\\
\emph{Assignment:}
\begin{equation*}
\p{P,p,T,t}=\texttt{get\_parameters\_function\_value\_suboptimality}(\texttt{algo},h,\alpha,j,\tau).
\end{equation*}
\emph{Arguments and allowed ranges:}
\begin{center}
\begin{minipage}{0.96\linewidth}
\centering
\texttt{algo}: \texttt{Algorithm};\quad
\(h,\alpha\in\naturals\);\quad
\(j\in\llbracket 1,\NumEval_{1}\rrbracket\);\quad
\(\tau\in\llbracket 0,h+\alpha+1\rrbracket\).
\end{minipage}
\end{center}
\emph{Induced targets:}
\begin{equation*}
\mathcal{V}\p{P,p,k}=0,
\qquad
\mathcal{R}\p{T,t,k}=f_{1}\p{y_{1,j}^{k+\tau}}-f_{1}\p{y^{\star}}.
\end{equation*}
}
\\
\bottomrule
\end{tabular}
\end{table}

We next present two examples showing how to use these interfaces and helper functions in practice.

\paragraph{Linear convergence: The Chambolle--Pock method.}
Consider the Chambolle--Pock method in \eqref{eq:chambolle_pock_interface:standard}, modeled in \texttt{AutoLyap} as \eqref{eq:chambolle_pock_interface:state_space}.
Assume that \(f_{1},f_{2}:\calH\to\reals\) are \(L\)-smooth and \(\mu\)-strongly convex. 
We seek the smallest certifiable \(\rho\in[0,1[\) such that
\begin{align}\label{eq:chambolle_pock:linear}
    \norm{x^k-y^{\star}}^{2} \in \mathcal{O}\p{\rho^{k}} \text{ as } k\to\infty,
\end{align}
where \(y^{\star}\in\zer\p{\partial f_1 + \partial f_2}\).
The code below illustrates this workflow for \((\mu,L,\tau,\sigma,\theta)=(0.05,50,1.6,1.6,0.22)\).
\begin{lstlisting}
from autolyap import IterationIndependent, SolverOptions
from autolyap.algorithms import ChambollePock
from autolyap.problemclass import InclusionProblem, SmoothStronglyConvex

mu, L = 0.05, 50.0
tau, sigma, theta = 1.6, 1.6, 0.22

problem = InclusionProblem([
    SmoothStronglyConvex(mu=mu, L=L),  # $f_1$
    SmoothStronglyConvex(mu=mu, L=L),  # $f_2$
])
algorithm = ChambollePock(tau=tau, sigma=sigma, theta=theta)

P, p, T, t = IterationIndependent.LinearConvergence.get_parameters_distance_to_solution(
    algorithm,
    i=1,
)

solver_options = SolverOptions(backend="mosek_fusion")
result = IterationIndependent.LinearConvergence.bisection_search_rho(
    problem,
    algorithm,
    P,
    T,
    p=p,
    t=t,
    solver_options=solver_options,
)

result["status"]  # value is "feasible"
result["rho"]     # approximately 0.8806
\end{lstlisting}
If \texttt{result["status"]} equals \texttt{"feasible"}, then \eqref{eq:chambolle_pock:linear} holds with \(\rho\) given by \texttt{result["rho"]}.
Repeating this over \(\tau=\sigma\in[1,1.8]\) and \(\theta\in[0,3/2]\) gives \Cref{fig:chambolle_pock:linear}.
\begin{figure}[!t]
\centering
  \begin{tikzpicture}
    \pgfplotstableread[col sep=comma]{data/linear_chambolle_pock/chambolle_pock_smooth_strongly_convex_tau_theta_rho.tex}\rhotable
    \begin{axis}[
      width=0.99\textwidth,
      height=0.50\textwidth,
      xlabel={\(\tau=\sigma\)},
      ylabel={\(\theta\)},
      ylabel style={rotate=-90},
      grid=both,
      xmin=0.99, xmax=1.75,
      ymin=0, ymax=1.5,
      scatter,
      scatter src=explicit,
      point meta min=0.88,
      point meta max=1.00,
      colorbar,
      colormap/viridis,
      colorbar horizontal,
      colorbar style={
          title={\(\rho\)},
          title style={
              at={(0,0)},
              anchor=east,
          },
      },
    ]
      \addplot [only marks, mark size=1.5pt] table [
        x=tau,
        y=theta,
        meta=rho,
        col sep=comma
      ] {\rhotable};
    \end{axis}
  \end{tikzpicture}
\caption{
  Smallest \(\rho\) certifiable via \texttt{AutoLyap} in \eqref{eq:chambolle_pock:linear} for the Chambolle--Pock method \eqref{eq:chambolle_pock_interface:standard}, where \(f_{1},f_{2}:\calH\to\reals\) are \(L\)-smooth and \(\mu\)-strongly convex with \((\mu,L)=(0.05,50)\).
  In the plot, we restrict to points such that \(\tau=\sigma \geq 1\).
  }
  \label{fig:chambolle_pock:linear}
\end{figure}

\paragraph{Sublinear convergence: The Chambolle--Pock method.}
Consider the Chambolle--Pock method in \eqref{eq:chambolle_pock_interface:standard}, modeled in \texttt{AutoLyap} as \eqref{eq:chambolle_pock_interface:state_space}.
Assume that \(f_{1},f_{2}:\calH\to\reals\) are proper, lower semicontinuous, and convex.
Suppose we want to verify whether the fixed-point residual \(\norm{\bx^{k+1}-\bx^{k}}^{2}\) is summable.
The code below verifies this for \(\p{\tau,\sigma,\theta,h,\alpha}=\p{1,1,1,1,1}\).
\begin{lstlisting}[label = lst:drs_mm_smlo]
from autolyap import IterationIndependent, SolverOptions
from autolyap.algorithms import ChambollePock
from autolyap.problemclass import Convex, InclusionProblem

tau, sigma, theta = 1.0, 1.0, 1.0
h, alpha = 1, 1

problem = InclusionProblem([
    Convex(),  # $f_1$
    Convex(),  # $f_2$
])
algorithm = ChambollePock(tau=tau, sigma=sigma, theta=theta)

P, p, T, t = IterationIndependent.SublinearConvergence.get_parameters_fixed_point_residual(
    algorithm, h=h, alpha=alpha
)

solver_options = SolverOptions(backend="mosek_fusion")
result = IterationIndependent.search_lyapunov(
    problem,
    algorithm,
    P,
    T,
    p=p,
    t=t,
    rho=1.0,
    h=h,
    alpha=alpha,
    solver_options=solver_options,
)

result["status"]  # value is "feasible"
result["certificate"]["Q"]  # Matrix $Q$
result["certificate"]["q"]  # Vector $q$ (None when m_func == 0)
result["certificate"]["S"]  # Matrix $S$
result["certificate"]["s"]  # Vector $s$ (None when m_func == 0)
\end{lstlisting}
If \texttt{result["status"]} equals \texttt{"feasible"}, then \(\p{\norm{\bx^{k+1}-\bx^{k}}^{2}}_{k\in\naturals}\) is summable.
Repeating this for multiple values of \(\tau=\sigma \in (1,2)\), \(\theta\in(0,3/2)\), and a few different values of \(\p{h,\alpha}\) gives \Cref{fig:chambolle_pock:sublinear}.
We observe that the region increases in size as \(h\) and \(\alpha\) increase.
However, in our numerical experiments, we did not observe any further increase beyond \((h, \alpha) = (2, 0)\), except for minor artefacts that may be attributed to the SDP solver.
Furthermore, the case \((h, \alpha) = (0, 0)\) corresponds to the smallest region, which matches the result in \cite[Figure 4a]{upadhyaya2025automatedtightlyapunov}, even though a different performance measure is used there.
\begin{figure}[!t]
\centering
  \begin{tikzpicture}
    \begin{axis}[
      width=0.99\textwidth,
      height=0.5\textwidth,
      xlabel={\(\tau=\sigma\)},
      ylabel={\(\theta\)},
      grid=both,
      xmax=1.75,
      xmin=0.99,
      ymin=0,
      ymax=1.5,
      ylabel style={rotate=-90},
    ]
      \pgfplotstableread[col sep=comma]{data/sublinear_chambolle_pock/chambolle_pock_fixed_point_residual_h2_alpha0.tex}\tabD
      \addplot [color=color5, only marks, mark size=1.5pt] table [x=x, y=y] {\tabD};
      \addlegendentry{\(\p{h,\alpha} = \p{2,0}\)}

      \pgfplotstableread[col sep=comma]{data/sublinear_chambolle_pock/chambolle_pock_fixed_point_residual_h1_alpha1.tex}\tabC
      \addplot [color=color3, only marks, mark size=1.5pt] table [x=x, y=y] {\tabC};
      \addlegendentry{\(\p{h,\alpha} = \p{1,1}\)}

      \pgfplotstableread[col sep=comma]{data/sublinear_chambolle_pock/chambolle_pock_fixed_point_residual_h1_alpha0.tex}\tabB
      \addplot [color=color2, only marks, mark size=1.5pt] table [x=x, y=y] {\tabB};
      \addlegendentry{\(\p{h,\alpha} = \p{1,0}\)}

      \pgfplotstableread[col sep=comma]{data/sublinear_chambolle_pock/chambolle_pock_fixed_point_residual_h0_alpha0.tex}\tabA
      \addplot [color=color1, only marks, mark size=1.5pt] table [x=x, y=y] {\tabA};
      \addlegendentry{\(\p{h,\alpha} = \p{0,0}\)}
    \end{axis}
  \end{tikzpicture}
\caption{
  Verification, via \texttt{AutoLyap}, of summability of the squared fixed-point residual \(\p{\norm{\bx^{k+1} - \bx^{k}}^{2}}_{k\in\naturals}\) for the Chambolle--Pock method \eqref{eq:chambolle_pock_interface:standard}, which is modeled by \eqref{eq:chambolle_pock_interface:state_space}. 
  Here, \(f_1,f_2:\calH\to\reals\cup\set{\pm\infty}\) are assumed to be proper, lower semicontinuous, and convex. 
  In the plot, we restrict to points such that \(\tau=\sigma \geq 1\), and we use a few different values of the history parameter \(h\) and the overlap parameter \(\alpha\).}
  \label{fig:chambolle_pock:sublinear}
\end{figure}

\subsection{Iteration-dependent Lyapunov analyses}\label{sec:iteration_dependent_Lyapunov}
In this section, we focus on finite-horizon analyses. These apply both to algorithms that run indefinitely and to algorithms with finite iteration budgets that match or exceed the horizon.
For notational simplicity, we assume in \eqref{eq:linear_system_with_nonlinearity} that \(K\in\mathbb{N}\) denotes the horizon.
Following \cite{taylor2019stochasticfirstorder}, we adopt an ansatz consisting of a sequence of iteration-dependent quadratic Lyapunov functions, also often referred to as \emph{potential functions} \cite{bansal2019potentialfunctionproofsgradientmethods}.

\begin{definition}\label{def:iteration_dependent_Lyapunov:chain}
Suppose that \(\bx_{0}\in\calH^{n}\) is an initial point, \(( ( \bx^{k},\allowbreak \bu^{k},\allowbreak \by^{k},\allowbreak \bFcn^{k} ) )_{k=0}^{K}\) is a sequence of iterates satisfying~\eqref{eq:linear_system_with_nonlinearity}, \(\p{y^{\star}, \allowbreak\hat{\bu}^{\star}, \allowbreak \bFcn^{\star}}\) is a point satisfying~\eqref{eq:solution}, \(Q_{k} \in \sym^{n+\NumEval+m}\) and \(q_{k} \in \reals^{\NumEvalFunc + \NumFunc}\) for each \(k\in\llbracket0,K\rrbracket\), and \(c_{K}\in\reals_{+}\).
Define
    \begin{equation}\label{eq:iteration_dependent_Lyapunov:chain:V}
        \p{\forall k \in \llbracket0,K\rrbracket} \quad \mathcal{V}(Q_{k},q_{k},k) = \quadform{Q_{k}}{\p{\bx^{k},\bu^{k},\hat{\bu}^{\star},y^{\star}}} + q_{k}^{\top}\p{\bFcn^{k},\bFcn^{\star}}.
    \end{equation}
We say that \((\p{Q_{k},q_{k}})_{k=0}^{K}\) and \(c_{K}\) satisfy a \emph{length-\(K\) sequence of chained Lyapunov inequalities} for algorithm \eqref{eq:linear_system_with_nonlinearity} over the problem class defined by \((\mathcal{F}_i)_{i\in\IndexFunc}\) and \((\mathcal{G}_i)_{i\in\IndexOp}\) if
    \begin{align}\label{eq:iteration_dependent_Lyapunov:chain}
        \mathcal{V}(Q_{K},q_{K},K) \leq \mathcal{V}(Q_{K-1},q_{K-1},K-1) \leq \ldots \leq \mathcal{V}(Q_{1},q_{1},1) \leq c_{K} \mathcal{V}(Q_{0},q_{0},0)
    \end{align}
holds for each initial point \(\bx_{0}\), for each sequence of iterates \(( ( \bx^{k},\allowbreak \bu^{k},\allowbreak \by^{k},\allowbreak \bFcn^{k} ) )_{k=0}^{K}\), for each point \(\p{y^{\star}, \hat{\bu}^{\star},\bFcn^{\star}}\), for each \(\p{f_{i}}_{i\in\IndexFunc} \in \prod_{i\in\IndexFunc} \mathcal{F}_{i}\), and for each \(\p{G_i}_{i\in\IndexOp} \in \prod_{i\in\IndexOp}\mathcal{G}_i\).
\end{definition}

In the proposed methodology, the user specifies \(\p{Q_{0},q_{0},Q_{K},q_{K}}\) and \texttt{AutoLyap} searches for intermediate \((\p{Q_{k},q_{k}})_{k=1}^{K-1}\) and the smallest \(c_{K}\) satisfying \eqref{eq:iteration_dependent_Lyapunov:chain}.
In particular, the user specifies the initial Lyapunov function \(\mathcal{V}(Q_{0},q_{0},0)\), final Lyapunov function \(\mathcal{V}(Q_{K},q_{K},K)\), and if \eqref{eq:iteration_dependent_Lyapunov:chain} holds, we can conclude that
\begin{align*}
    \mathcal{V}(Q_{K},q_{K},K) \leq c_{K} \mathcal{V}(Q_{0},q_{0},0).
\end{align*}

\subsubsection{Interface and example}

In \texttt{AutoLyap}, the main interface for \Cref{def:iteration_dependent_Lyapunov:chain} is \texttt{IterationDependent.search\_lyapunov}.
For fixed \(\p{Q_{0},q_{0},Q_{K},q_{K},K}\), it formulates the SDP associated with \Cref{def:iteration_dependent_Lyapunov:chain}, with the objective of minimizing \(c_{K}\), and solves it using the backend selected through \texttt{SolverOptions}.

To support this workflow, \texttt{AutoLyap} provides helper functions in \path{IterationDependent} for choosing endpoint parameters \(\p{Q_{0},q_{0}}\) and \(\p{Q_{K},q_{K}}\).
\Cref{tab:iter_dep_helpers} summarizes some built-in \path{get_parameters_*} helper functions in \path{IterationDependent}.

\begin{table}[t]
\centering
\caption{
Some helper functions in \texttt{IterationDependent}, with inputs, allowed values, and induced targets \(\mathcal{V}\p{Q_{k},q_{k},k}\).
Notation follows \eqref{eq:linear_system_with_nonlinearity} and \Cref{def:iteration_dependent_Lyapunov:chain}, with \(\mathcal{V}\) defined in \eqref{eq:iteration_dependent_Lyapunov:chain:V}.
When \(\NumFunc = 0\), the vector \(q_{k}\) is omitted (equivalently, helpers return \(Q_{k}\) instead of \(\p{Q_{k},q_{k}}\)).
Further details are found in the \href{https://autolyap.github.io}{documentation}.
}
\label{tab:iter_dep_helpers}
\begin{tabular}{p{0.93\textwidth}}
\toprule
\rowcolor{black!8}[\tabcolsep][\tabcolsep]
\parbox{\linewidth}{
\emph{Assignment:}
\begin{equation*}
\p{Q_{k},q_{k}}=\texttt{get\_parameters\_distance\_to\_solution}(\texttt{algo},k,i,j).
\end{equation*}
\emph{Arguments and allowed ranges:}
\begin{center}
\begin{minipage}{0.96\linewidth}
\centering
\texttt{algo}: \texttt{Algorithm};\quad
\(k\in\llbracket 0,K\rrbracket\);\quad
\(i\in\llbracket 1,m\rrbracket\);\quad
\(j\in\llbracket 1,\NumEval_{i}\rrbracket\).
\end{minipage}
\end{center}
\emph{Induced target:}
\begin{equation*}
\mathcal{V}\p{Q_{k},q_{k},k}=\norm{y_{i,j}^{k}-y^{\star}}^{2}.
\end{equation*}
}
\\
\addlinespace[0.6em]
\rowcolor{white}[\tabcolsep][\tabcolsep]
\parbox{\linewidth}{
\emph{Requires:} \(m=\NumFunc=1\).\\
\emph{Assignment:}
\begin{equation*}
\p{Q_{k},q_{k}}=\texttt{get\_parameters\_function\_value\_suboptimality}(\texttt{algo},k,j).
\end{equation*}
\emph{Arguments and allowed ranges:}
\begin{center}
\begin{minipage}{0.96\linewidth}
\centering
\texttt{algo}: \texttt{Algorithm};\quad
\(k\in\llbracket 0,K\rrbracket\);\quad
\(j\in\llbracket 1,\NumEval_{1}\rrbracket\).
\end{minipage}
\end{center}
\emph{Induced target:}
\begin{equation*}
\mathcal{V}\p{Q_{k},q_{k},k}=f_{1}\p{y_{1,j}^{k}}-f_{1}\p{y^{\star}}.
\end{equation*}
}
\\
\addlinespace[0.6em]
\rowcolor{black!8}[\tabcolsep][\tabcolsep]
\parbox{\linewidth}{
\emph{Assignment:}
\begin{equation*}
\p{Q_{k},q_{k}}=\texttt{get\_parameters\_fixed\_point\_residual}(\texttt{algo},k).
\end{equation*}
\emph{Arguments and allowed ranges:}
\begin{center}
\begin{minipage}{0.96\linewidth}
\centering
\texttt{algo}: \texttt{Algorithm};\quad
\(k\in\llbracket 0,K\rrbracket\).
\end{minipage}
\end{center}
\emph{Induced target:}
\begin{equation*}
\mathcal{V}\p{Q_{k},q_{k},k}=\norm{\bx^{k+1}-\bx^{k}}^{2}.
\end{equation*}
}
\\
\bottomrule
\end{tabular}
\end{table}

\begin{remark}
    A key computational feature is that each inequality in \eqref{eq:iteration_dependent_Lyapunov:chain} is enforced through a one-step SDP condition of fixed size.
    Consequently, as the horizon \(K\) increases, the total number of constraints and decision variables grows linearly in \(K\).
\end{remark}

\paragraph{Nesterov's fast gradient method.}
Consider Nesterov's fast gradient method given by \eqref{eq:nesterov_fast_interface:standard_form}, which is modeled as \eqref{eq:nesterov_fast_interface:state_space} in \texttt{AutoLyap}.
Suppose we want to certify a finite-horizon bound of the form
\begin{equation}\label{eq:nesterov_fast_iteration_dependent_bound}
f_{1}\p{x^{K}} - f_{1}\p{y^{\star}} \le c_{K}\norm{x^{0}-y^{\star}}^{2},
\end{equation}
where \(y^{\star}\in\zer\p{\partial f_1}\), for some horizon \(K\in\mathbb{N}\).
The code below shows this workflow for \(\p{L,\gamma,K}=\p{1,1,10}\).
\begin{lstlisting}
from autolyap import IterationDependent, SolverOptions
from autolyap.algorithms import NesterovFastGradientMethod
from autolyap.problemclass import InclusionProblem, SmoothConvex

L, gamma, K = 1.0, 1.0, 10

problem = InclusionProblem([SmoothConvex(L)])
algorithm = NesterovFastGradientMethod(gamma=gamma)

Q_0, q_0 = IterationDependent.get_parameters_distance_to_solution(
    algorithm, k=0, i=1, j=2
)
Q_K, q_K = IterationDependent.get_parameters_function_value_suboptimality(
    algorithm, k=K, j=2
)

solver_options = SolverOptions(backend="mosek_fusion")
result = IterationDependent.search_lyapunov(
    problem,
    algorithm,
    K,
    Q_0,
    Q_K,
    q_0=q_0,
    q_K=q_K,
    solver_options=solver_options,
)

result["status"]  # value is "feasible"
result["c_K"]     # approximately $0.01103$
k = 3 # intermediate iteration index
result["certificate"]["Q_sequence"][k]  # $Q_k$ for k = 3
result["certificate"]["q_sequence"][k]  # $q_k$ for k = 3
\end{lstlisting}
Repeating this for \(K \in \llbracket1,100\rrbracket\) gives \Cref{fig:nesterov_fast}.
\begin{figure}[t]
\centering
  \begin{tikzpicture}
    \begin{loglogaxis}[
      width=0.99\textwidth,
      height=0.5\textwidth,
      grid=both,
      grid style={line width=.1pt, draw=gray!30},
      major grid style={line width=.2pt, draw=gray!50},
      xlabel={$K$},
      ylabel={$c_K$},
      ylabel style={rotate=-90},
      legend pos=north east,      
      legend cell align=left      
    ]

      \addplot+[
        color=color1,
        solid,
        mark=none,
        line width=1.2pt,
      ]
      table [x=x, y=y, col sep=comma]
      {data/nesterov_fast_gradient_method/automatic_lyapunov.tex};
      \addlegendentry{\texttt{AutoLyap}}

      \addplot+[
        color=color8,
        dash pattern=on 4pt off 2pt,
        mark=none,
        line width=1.2pt,
      ]
      table [x=x, y=y, col sep=comma]
      {data/nesterov_fast_gradient_method/nesterov_second_bound.tex};
      \addlegendentry{$L/(2\lambda_K^2)$}

      \addplot+[
        color=color5,
        dashed,
        mark=none,
        line width=1.2pt,
      ]
      table [x=x, y=y, col sep=comma]
      {data/nesterov_fast_gradient_method/nesterov_first_bound.tex};
      \addlegendentry{$2L/(K+2)^2$}

    \end{loglogaxis}
  \end{tikzpicture}
\caption{Constants \(c_{K}\) in \eqref{eq:nesterov_fast_iteration_dependent_bound} for Nesterov's fast gradient method \eqref{eq:nesterov_fast_interface:standard_form} with step size \(\gamma=1\), applied to a convex and \(L\)-smooth function \(f_{1}:\calH\to\reals\) with \(L=1\). Values are obtained with \texttt{AutoLyap} and compared with the classical rates in \cite{nesterov1983fast}.}
  \label{fig:nesterov_fast}
\end{figure}

\section{The \texttt{Julia} package \texttt{AutoLyap.jl}}\label{sec:autolyap_jl}

This section presents the package \texttt{AutoLyap.jl}, a pure-\texttt{Julia}, open-source implementation of the methodology described in this paper, available at 
\begin{center}
    \url{https://github.com/autolyap/autolyap.jl}.
\end{center}
\texttt{AutoLyap.jl} matches the \texttt{Python} package \texttt{AutoLyap} feature-for-feature and has nearly identical syntax, making it easy for end users to switch seamlessly between the two packages.
For example, the following \texttt{Julia} code illustrates the linear convergence rates of the Douglas–Rachford method for the inclusion problem \(0 \in G_1(y) + G_2(y)\), where \(G_1: \mathcal{H} \rightrightarrows \mathcal{H}\) is maximally monotone and \(G_2: \mathcal{H} \to \mathcal{H}\) is \(\mu\)-strongly monotone and \(L\)-Lipschitz continuous, with \(\mu = 1\) and \(L = 2\).

\begin{lstlisting}[style=juliastyle]
using AutoLyap
using AutoLyap: IterationIndependent

components_list = [MaximallyMonotone(), # $G_1$
    [
	StronglyMonotone(mu = 1.0),
	LipschitzOperator(L = 2.0)] # $G_2$
]
problem = InclusionProblem(components_list)
algorithm = DouglasRachford(gamma = 1.0, lambda_value = 2.0, operator_version=true)
(P, T) =  IterationIndependent.get_parameters_distance_to_solution(algorithm)
result = IterationIndependent.bisection_search_rho(problem, algorithm, P, T)
rho = result["rho"] # value of $\rho$ is approximately 0.4286
\end{lstlisting}

While the two packages are functionally equivalent, they differ in their architectural design.

\paragraph{Object-oriented programming in \texttt{Python} vs.\ multiple dispatch in \texttt{Julia}.}

The core architectural difference lies in the programming paradigm each package employs.
The \texttt{Python} package uses \emph{object-oriented programming}, whereas the \texttt{Julia} package utilizes the \emph{multiple dispatch} paradigm.
Neither paradigm is inherently superior: \texttt{Python} benefits from a broader ecosystem, whereas \texttt{Julia} is designed for high-performance scientific computing.

Object-oriented programming in \texttt{Python} revolves around the concept of an \emph{object} that encapsulates both data and methods within a \emph{class}.
\texttt{Python} enables (i) \emph{inheritance} by allowing subclasses to extend or override parent class functionality and (ii) \emph{polymorphism}, i.e., the ability of the same method name to have multiple implementations through method overriding.
As a result, method dispatch in \texttt{Python} is effectively single-dispatch: it depends primarily on the type of the receiver, \texttt{self}.

In contrast, the \emph{multiple dispatch} paradigm used by \texttt{Julia} decouples data from methods.
The data structure called \texttt{struct} bundles data-related fields and can be used to create composite types, but does not contain any methods other than inner constructors.
The methods live outside of the \texttt{struct}s in functions (a definition of one possible behavior for a function is called a method in \texttt{Julia}) and can admit multiple definitions based on the types of all arguments; \texttt{Julia} selects the most specific applicable method at runtime.
For code reuse and hierarchy, \texttt{Julia} intentionally avoids the notion of inheritance, and achieves the same functionality through (i) \emph{composition}, where instead of inheriting from a type, an instance of that type is included in a new \texttt{struct}, and (ii) an \emph{abstract type hierarchy}, in which a concrete \texttt{struct} can be a subtype of an abstract type.
Multiple dispatch in \texttt{Julia} provides \emph{a high degree of polymorphism}: any function can be extended with new methods for any data type (including types from other libraries) without requiring changes to existing code.

\paragraph{Modeling language for the SDPs.}



\texttt{AutoLyap.jl} is built on \texttt{JuMP} \cite{lubin2023jumprecentimprovementsmodelinglanguage}---an open-source, solver-agnostic modeling language for mathematical optimization embedded in \texttt{Julia}, with syntax that enables the natural algebraic representation of optimization problems and deep integration with the high-performance \texttt{Julia} ecosystem.
Hence, a model defined in \texttt{AutoLyap.jl} is not tied to any specific solver, and the user can switch from a commercial solver like \texttt{MOSEK} to an open-source solver like \texttt{Clarabel} \cite{clarabel2024interiorpointsolverconicprograms} by changing a single line of code in \texttt{AutoLyap.jl}:

\begin{lstlisting}[style=juliastyle]
solver = :clarabel 
\end{lstlisting}

At present, \texttt{AutoLyap.jl} can access a wide range of state-of-the-art SDP solvers through \texttt{JuMP}, including \texttt{Clarabel} \cite{clarabel2024interiorpointsolverconicprograms} (open-source, \texttt{AutoLyap.jl} uses it by default), \texttt{MOSEK} (commercial, free for academic use), \texttt{COPT} \cite{copt2023cardinaloptimizeruserguide} (commercial, free for academic use), \texttt{SCS} \cite{scs2016conicoptimizationoperatorsplittinghomogeneousselfdualembedding} (open-source), \texttt{COSMO} \cite{garstka2021cosmoconicoperatorsplittingmethod} (open-source), \texttt{SDPA} \cite{YamashitaEtal2012_SDPAFamily} (open-source), \texttt{ProxSDP} \cite{souto2020exploiting} (open-source), \texttt{Hypatia} \cite{coey2022solving} (open-source), and \texttt{SDPLR} \cite{burer2006computational} (open-source).
Hence, \texttt{AutoLyap.jl} can be a suitable choice in academic and research-oriented environments with open-source development goals or for users who lack licenses for commercial solvers.

\section{Conclusion}
\label{sec:conclusions}

This paper introduced AutoLyap, a software suite designed to assist with Lyapunov analyses of first-order methods via semidefinite programming.
It simplifies the derivation and verification of convergence guarantees for structured optimization and inclusion problems.

AutoLyap is designed to encourage community involvement and to enable researchers to extend the framework collaboratively.
We mention two relevant directions for future development.
\begin{enumerate}[(i)]
\item The algorithm representation~\eqref{eq:linear_system_with_nonlinearity} can be extended to allow for more types of oracles (including, e.g., Frank--Wolfe-type oracles~\cite{taylor2017exactworstcase}, Bregman-type oracles~\cite{dragomir2021optimalcomplexitycertification}, or approximate oracles with explicit error bounds~\cite{barre2022principledanalysesdesign}).
\item Ideally, one would like to move beyond relying on numerical results as an intermediate step and instead pursue closed-form Lyapunov analyses directly. In our context, this amounts to finding closed-form solutions to parametric semidefinite programs, which is an active area of research \cite{basu2006algorithmsrealalgebraic,naldi2025solvinggenericparametric}.
\end{enumerate}

\section*{Acknowledgments}
M. Upadhyaya and A.B. Taylor are supported by the European Union (ERC grant CASPER 101162889).
The French government also partly supported this work through the Agence Nationale de la Recherche as part of the ``France 2030'' program, reference ANR-23-IACL-0008 ``PR[AI]RIE-PSAI''.
M. Upadhyaya and P. Giselsson acknowledge support from the ELLIIT Strategic Research Area and the Wallenberg AI, Autonomous Systems, and Software Program (WASP), funded by the Knut and Alice Wallenberg Foundation.
Additionally, P. Giselsson acknowledges support from the Swedish Research Council. S.\ Das Gupta acknowledges support from AFOSR Grant Number FA9550-25-1-0183.
Views and opinions expressed are, however, those of the authors only and do not necessarily reflect those of the funding agencies or the granting authorities, who cannot be held responsible for them.

\begingroup
\sloppy
\printbibliography
\endgroup

\end{document}